\documentclass{amsart}
\usepackage[english]{babel}

\newtheorem{thm}{Theorem}[subsection]
\newtheorem{prop}[thm]{Proposition}
\newtheorem{lem}[thm]{Lemma}
\newtheorem{cor}[thm]{Corollary}
\theoremstyle{definition}
\newtheorem{dfn}[thm]{Definition}
\theoremstyle{remark}
\newtheorem{rmk}[thm]{Remark}

\numberwithin{equation}{subsection}

\newcommand{\x}{\times}
\newcommand{\ox}{\otimes}
\renewcommand{\tt}{{\gt}}  
\newcommand{\abs}[1]{\lvert#1\rvert}
\newcommand{\norm}[1]{\lVert#1\rVert}
\def\lacute{\mathopen{<}}\def\racute{\mathopen{>}}
\newcommand{\scal}[2]{{\lacute#1,#2\racute}}
\newcommand{\ensemble}[2]{\{\,#1\mid#2\,\}}

\newcommand{\C}{{\mathbb C}}
\renewcommand{\L}{{\mathbb L}}
\newcommand{\N}{{\mathbb N}}

\newcommand{\R}{{\mathbb R}}
\newcommand{\Z}{{\mathbb Z}}

\newcommand{\CaD}{\mathcal D}   
       
\newcommand{\CF}{\mathcal F}

\newcommand{\CI}{\mathcal I}
\newcommand{\CJ}{\mathcal J}

\newcommand{\CM}{\mathcal M}
\newcommand{\CN}{\mathcal N}

\newcommand{\CO}{\mathcal O}

\newcommand{\kN}{\mathfrak N}
\newcommand{\kO}{\mathfrak O}

\newcommand{\ka}{\mathfrak a} \newcommand{\tka}{\tilde\ka}

\newcommand{\kc}{\mathfrak c}
\newcommand{\kg}{\mathfrak g} \newcommand{\tkg}{\tilde\kg}
\newcommand{\kh}{\mathfrak h}
\newcommand{\kk}{\mathfrak k} \newcommand{\tkk}{\tilde\kk}
\newcommand{\kl}{\mathfrak l}
\newcommand{\km}{\mathfrak m}

\newcommand{\kp}{\mathfrak p} \newcommand{\tkp}{\tilde\kp}
\newcommand{\kq}{\mathfrak q}
\newcommand{\ks}{\mathfrak s}
\newcommand{\Dp}{\CaD_\kp}
\newcommand{\Dps}{\CaD_{\kp^s}}

\newcommand{\kgp}{\kg_P}
\newcommand{\kqp}{\kq_P}
\newcommand{\khp}{\kh_P}
\newcommand{\khpo}{\khp^\bot}
\newcommand{\khpp}{(\khpo)'}
\newcommand{\kap}{\ka_P}
\newcommand{\kapo}{\kap^\bot}
\newcommand{\kapd}{\ka_{P_2}}
\newcommand{\kapdo}{\kapd^\bot}
\newcommand{\kapp}{(\kapo)'}
\newcommand{\sld}{\ks\kl_2}

\newcommand{\lp}{\gl_\kp}
\newcommand{\lpx}{\gl_\kp(x)}

\newcommand{\arr}{\xrightarrow{\ \ }}
\renewcommand{\setminus}{-}     

\renewcommand{\Re}{\operatorname{Re}}

\newcommand{\dsur}[1]{\frac \partial{\partial#1} }               
\newcommand{\TX}{{T^*X}}

\newcommand{\tYX}{{T_Y\!X}}

\newcommand{\TYX}{{T^*_YX}}

\newcommand{\OX}{{{\mathcal O_X}}}
\newcommand{\OY}{{\mathcal O_Y}}

\newcommand{\DX}{{\mathcal D_X}}
\newcommand{\DXi}[1]{{\mathcal D_{X,#1}}}

\newcommand{\DY}{{{\mathcal D_Y}}}
\newcommand{\DZ}{{{\mathcal D_Z}}}
\newcommand{\DYX}{{\mathcal D_{Y\rightarrow X}}}
\newcommand{\DXY}{{\mathcal D_{X\leftarrow Y}}}

\newcommand{\DtYXa}{{\mathcal D_{[\tYX]}}}
\newcommand{\DtYXr}{{\mathcal D_{[\tYX/Y]}}}

\newcommand{\DIXY}{{\mathcal D^\infty_{X\leftarrow Y}}}
\newcommand{\DIY}{{{\mathcal D^\infty_Y}}}

\newcommand{\tL}{\widetilde L}
\newcommand{\tM}{\widetilde M}

\newcommand{\tX}{\widetilde X}
\newcommand{\tY}{\widetilde Y}
\newcommand{\tZ}{\widetilde Z}

\newcommand{\CHM}{{Ch(\CM)}}

\newcommand{\DbM}{{\mathcal D\!b_M}}
\newcommand{\DbN}{{\mathcal D\!b_N}}
\newcommand{\gNDbM}{{\Gamma_N(\DbM)}}

\DeclareMathOperator{\Ext}{\mathcal Ext}

\DeclareMathOperator{\HOM}{\mathcal Hom}
\DeclareMathOperator{\End}{\mathcal End}
\DeclareMathOperator{\RHOM}{\R\mathcal Hom}

\DeclareMathOperator{\redim}{dim_0}
\DeclareMathOperator{\rank}{rank}
\DeclareMathOperator{\Ker}{Ker}

\DeclareMathOperator{\ad}{ad}

\newcommand{\ga}{\alpha}
\newcommand{\gb}{\beta}
\newcommand{\gga}{\gamma}           \newcommand{\gG}{\Gamma}
\newcommand{\gd}{\delta}           \newcommand{\gD}{\Delta}
                                   \newcommand{\tgD}{\widetilde\Delta}
\newcommand{\gep}{\varepsilon}

\newcommand{\gh}{\eta}
\newcommand{\gth}{\theta}          \newcommand{\gTh}{\Theta}
\newcommand{\gvt}{\vartheta}
\newcommand{\gk}{\kappa}
\newcommand{\gl}{\lambda}          
\newcommand{\gm}{\mu}

\newcommand{\gx}{\xi}              
\newcommand{\gp}{\pi}              

\newcommand{\gro}{\varrho}    
\newcommand{\gs}{\sigma}           \newcommand{\gS}{\Sigma}
\newcommand{\gt}{\tau}
\newcommand{\gf}{\varphi}          \newcommand{\gF}{\Phi}
\newcommand{\gq}{\chi}
             
           \newcommand{\gO}{\Omega}

\newcommand{\tgh}{{\widetilde\eta}}

\newcommand{\Dg}{\CaD_\kg}

\newcommand{\MF}{\CM_F}
\newcommand{\MFf}{\widehat{\CM_F}}

\newcommand{\MFl}{\CM_\gl^{\mathcal F}}

\newcommand{\kgrs}{\kg_{rs}}
\newcommand{\kprs}{\kp_{rs}}

\begin{document}

\title[Characters of semisimple groups]{$\mathcal D$-Modules and Characters of Semisimple Lie Groups}
\author{Esther Galina}
\address{CIEM-FAMAF\\
    Universidad Nacional de C\'ordoba\\
    Ciudad Universitaria\\
    5000 C\'ordoba\\
    ARGENTINA}
\email{galina@mate.uncor.edu}
\author{Yves Laurent}
\address{ Institut Fourier Math\'ematiques \\
    UMR 5582  CNRS/UJF\\
       BP 74\\
F-38402 St Martin d'H\`eres Cedex\\
  FRANCE }
\email{Yves.Laurent@ujf-grenoble.fr}
\urladdr{http://www-fourier.ujf-grenoble.fr/\~{ }laurenty}
\thanks{supported by CONICET,CONICOR, FONCYT, TWAS, UNC and University
  of Grenoble}

\begin{abstract}
A famous theorem of Harish-Chandra asserts that all invariant
eigendistributions on a semisimple Lie group are locally
integrable functions. We show that this result and its extension
to symmetric pairs are consequences of an algebraic property of a
holonomic $\mathcal D$-module defined by Hotta and Kashiwara.
\end{abstract}

\maketitle

\tableofcontents

\section*{Introduction}

Let $G_\R$ be a real semisimple Lie group. An \textsl{invariant
eigendistribution} on $G_\R$ is a distribution which is invariant
under the adjoint action of $G_\R$ and an eigenvalue of every
biinvariant differential operator on $G_\R$. Any irreducible
representation of $G_\R$ has a character which is an invariant
eigendistribution.

A famous theorem of Harish-Chandra \cite{HC1} asserts that all invariant
eigendistributions are locally integrable functions. The classical proof of
the theorem is divided in three steps:

(i)Any invariant eigendistribution $\gq$ is analytic on the set $G_{rs}$ of
regular semisimple points which is a Zariski dense open subset of $G_\R$.

(ii)The restriction $F$ of $\gq$ to $G_{rs}$ extends to a $L^1_{loc}$ function
on $G_\R$.

(iii)There is no invariant eigendistribution supported by $G_\R\setminus
G_{rs}$. \vskip 0.2cm

The problem is local and invariant eigendistributions may be studied
on the Lie algebra $\kg_\R$ of $G_\R$. Hotta and Kashiwara defined in
\cite{HOTTA} a family of holonomic $\CaD$-modules on the
complexification $\kg$ of $\kg_\R$ whose solutions are the invariant
eigendistributions. This module is elliptic on the set $\kg_{rs}$ of
regular semisimple points which shows (i) and using the results of
Harish-Chandra they also proved that it is regular.

In \cite{SEKI}, J. Sekiguchi extended partially these results to
symmetric pairs. There is an analog to the $\CaD$-module of
Hotta-Kashiwara, which is holonomic and elliptic on the regular
semisimple points. But the result of Harish-Chandra do not always
extend and Sekiguchi gave a counter-example. He introduced a class of
symmetric pairs ("nice pairs") for which he proved (iii), that is no
distribution solution is supported by $\kg_\R\setminus\kg_{rs}$. He
also extended the result to hyperfunctions and conjectured that the
Hotta-Kashiwara $\CaD$-module is regular in the case of symmetric
pairs. In \cite{LBY}, we proved this conjecture for all symmetric
pairs. This shows, among others, that all hyperfunction solutions are
distributions.

In several papers \cite{LEVASS}\cite{LEVASS3}\cite{LEVASS2}, Levasseur and
Stafford gave an algebraic proof of point (iii) for distributions in the case
of semisimple groups and in the case of nice symmetric pairs.

The aim of this paper is to show that the Harish-Chandra theorem and
its extension to symmetric pairs is a consequence of an algebraic
property of the Hotta-Kashiwara $\CaD$-module.  This property is the
following:

If $\CM$ is a holonomic $\DX$-module on a manifold $X$, to each
submanifold $Y$ of $X$ is associated a polynomial which is called the
$b$-function of $\CM$ along $Y$ (see \S\ref{sec:bfunct} for a precise
definition). We say that the module $\CM$ is \textsl{tame} if there
exists a locally finite stratification $X=\cup X_\ga$ such that, for
each $\ga$, the roots of the $b$-function of $\CM$ along $X_\ga$ are
greater than the opposite of the codimension of $X_\ga$.

We show first that the distribution solutions of a tame $\CaD$-module satisfy
properties (i)-(ii)-(iii) (replacing $G\setminus G_{rs}$ by the singular
support of $\CM$) and second that the Hotta-Kashiwara module is tame. More
precisely, we show that it is always tame in the semi-simple case and in the
case of symmetric pairs, we find a relation between the roots of the
$b$-functions and some numbers introduced by Sekiguchi. This relation implies
that the module is tame for nice pairs. In fact, this is true after an
extension of the definition of  tame $\CaD$-module which is given in
section \ref{sec:tameres}.

In this way, we get a new proof of the results of Harish-Chandra,
Sekiguchi and Levasseur-Stafford.
Concerning the integrability of solutions in the case of symmetric
pairs, no result was known. We get this integrability but we need a
condition which is slightly stronger than the condition satisfied by
nice pairs (see section \ref{sec:Sym}).

Tame $\CaD$-modules have other nice properties, in particular they
have no quotients supported by a hypersurface. In the complex domain,
a Nilsson class solution is always a $L^2_{loc}$-function.

From the point of view of $\CaD$-modules, our work establish a new
kind of connection between algebraic properties of a holonomic
$\CaD$-module and the growth of its solutions. A result of
Kashiwara \cite{KASHREG} shows that hyperfunction solutions of a
regular $\CaD$-module are distributions. On the other hand,
results of Ramis \cite{RAMIS} in the one dimensional case and of
one of the authors \cite{IRR} in the general case, rely on the
Gevrey or exponential type of the solutions with the Newton
Polygon of the $\CaD$-module. Here we show that the $L^p$ growth
of the solutions is given by the roots of the $b$-functions.

Another interest of our work is to give an example of a family of non trivial
holonomic $\CaD$-modules for which it is possible to calculate
explicitly the $b$-functions.

This paper is divided in three parts. The first section is devoted to
the definitions and the statement of the main results. We recall the
definition of $b$-functions in section \ref{sec:bfunct}, but this is
not sufficient in the case of the Hotta-Kashiwara module and we have
to extend slightly this definition in section \ref{sec:gbfunct}. Then
we define the ``tame''-$\CaD$-modules and give our principal
results.

In the second part, we study the relations between the roots
of the $b$-functions and the growth of solutions, proving in
particular that the distribution solutions of a tame $\CaD$-module are
locally integrable.

In the third section, we calculate the $b$-functions of the Hotta-Kashiwara
module and show that it is tame. The key point of the proof is the fact that
the Fourier transform of the Hotta-Kashiwara module is supported by the
nilpotent cone.

\section{The main results}  \label{chap1}
\subsection{$\CaD$-modules and generators}\label{sec:generators}

Let $(X,\CO_{[X]})$ be a smooth algebraic variety defined over $\C$ and
$(X,\OX)$ be the underlying complex manifold. We will denote by
$\CaD_{[X]}$ the sheaf of differential operators with coefficients in
$\CO_{[X]}$ and $\DX$ be the sheaf of differential operators with
coefficients in $\OX$. The theories of $\CaD_{[X]}$-modules and
$\DX$-modules are very similar, and we refer to \cite{HOTTA} for an
introduction to holonomic and regular holonomic $\CaD_{[X]}$-modules.
In this paper, we will work mostly with $\DX$-modules but the
definitions of section \ref{chap1} are valid in both cases.

In the theory of $\DX$-modules, what is called the ``solutions of
a coherent $\DX$-module $\CM$ in a sheaf of functions $\CF$'' is
the derived functor $\RHOM_\DX(\CM,\CF)$. In this paper, we will
be interested only in its first cohomology group, that is the
sheaf $\HOM_\DX(\CM,\CF)$.  If $\CM$ is a cyclic $\DX$-module, the
choice of a generator defines an isomorphism $\CM\simeq\DX/\CI$
where $\CI$ is a coherent ideal of $\DX$. Then there is a
canonical isomorphism between $\HOM_\DX(\CM,\CF)$ and
$\ensemble{u\in\CF}{\forall P\in\CI,\ \ Pu=0}$.  But this
isomorphism depends on the choice of a generator of $\CM$ and some
properties of solutions of partial differential equations as to be
a $L^2$-function depend on this choice. So, we will always
consider $\DX$-modules explicitly written as $\DX/\CI$ or
$\DX$-modules with a given set of generators for which there is no
ambiguity.

A similar situation will be found when defining the $b$-functions
in section \ref{sec:bfunct}.

\subsection{$V$-filtration}\label{sec:vfiltr}

Let $Y$ be a smooth subvariety of $X$.  The sheaf $\DX$ is provided
with two canonical filtrations. First, we have the filtration by the usual
order of operators denoted by $(\DXi m)_{m\ge0}$ and second the
$V$-filtration of Kashiwara \cite{KVAN}~:
$$V_k\DX=\ensemble{P\in\DX}{\forall j\in\Z,\ \ P\CJ_Y^j\subset\CJ_Y^{j-k}}$$
where $\CJ_Y$ is the definition ideal of $Y$ and $\CJ_Y^j=\OX$ if
$j\le0$.

In coordinates $(x,t)$ such that $Y=\{t=0\}$, $\CJ_Y^k$ is, for
$k\ge0$, the sheaf of functions
$$\sum_{\abs{\ga}=k}f_\ga(x,t)t^\ga$$
hence the operators $x_i$ and $D_{x_i}:=\dsur{x_i}$ have order $0$ for
the $V$-filtration while the operators $t_i$ have order $-1$ and $D_{t_i}:=\dsur{t_i}$ order
$+1$.

The associated graduate is defined as~:
$$gr_V\DX:=\oplus gr^k_V\DX,\qquad
gr^k_V\DX:=V_k\DX\left/V_{k-1}\DX\right.$$
By definition,
$gr_V\DX$ operates on the direct sum $\bigoplus
\left(\CJ_Y^k\left/\CJ_Y^{(k+1)}\right.\right)$. But this sheaf is
canonically isomorphic to the direct image by the projection
$p:\tYX\to Y$ of the sheaf  $\CO_{[\tYX]}$ of holomorphic
functions on the normal bundle $\tYX$ polynomial in the fibers of
$p$ (in the algebraic case, it is the sheaf $\CO_{\tYX}$ of
functions on $\tYX$). In this way $gr_V\DX$ is a subsheaf of
$p_*\HOM_\C(\CO_{[\tYX]},\CO_{[\tYX}])$ and it is easily verified
in coordinates that this subsheaf is exactly the sheaf of
differential operators with coefficients in $\CO_{[\tYX]}$:
$$gr_V\DX\simeq p_*\DtYXa$$

The graduate associated to the filtration $(\DXi m)$
is $gr\DX\simeq \gp_*\CO_{T^*X}$ where $\gp:T^*X\to X$ is the
cotangent bundle in the algebraic case and the sheaf
$\gp_*\CO_{[T^*X]}$ of holomorphic functions polynomial in the fibers
of $\gp$ in the analytic case.

Let $\CM$ be a coherent $\DX$-module. A good filtration of $\CM$ is a
filtration which is locally finitely generated that is locally of the
form~:
$$\CM_m=\sum_{j=1,\dots,N}\DXi {m+m_j}u_j$$
where $u_1,\dots,u_N$ are (local) sections of $\CM$ and
$m_1,\dots,m_N$ integers.

It is well known that if $(\CM_m)$ is a good filtration of $\CM$,
the associated graduate $gr\CM$ is a coherent $gr\DX$-module,
that is a coherent $\gp_*\CO_{[T^*X]}$-module and defines the
characteristic variety of $\CM$ which is a subvariety of $T^*X$.
This subvariety is involutive for the canonical symplectic
structure of $T^*X$ and a $\DX$-module is said to be holonomic if
its characteristic variety is lagrangian that is of minimal
dimension. If $\CM$ is holonomic, its sheaf of endomorphisms
$\End_\DX(\CM)=\HOM_\DX(\CM,\CM)$ is constructible, that is there
is a stratification of $X$ for which $\End_\DX(\CM)$ is locally
constant and finite dimensional on each stratum.

In the same way, a good $V$-filtration of $\CM$ is a filtration of
$\CM$ associated to the $V$-filtration of $\DX$ which is locally
finitely generated that is locally~:
$$V_k\CM=\sum_{j=1,\dots,N}(V_{k+k_j}\DX) u_j$$

If $V\CM$ is a good filtration, the associated graduate
$gr_V\CM$ is a coherent $gr_V\DX$-module hence $p^{-1}gr_V\CM$ is
a coherent $\DtYXa$-module \cite{KVAN}. Moreover it may be proved
that if $\CM$ is holonomic then $p^{-1}gr_V\CM$ is a holonomic
$\DtYXa$-module \cite{ENS} hence $\End_\DtYXa(p^{-1}gr_V\CM)$ is locally
finite dimensional.

\subsection{$b$-functions}\label{sec:bfunct}

The fiber bundle $\tYX$ is provided with a canonical vector field, the
Euler vector field $\gvt$ characterized by $\gvt f=kf$ for any
function $f$ homogeneous of degree $k$ in the fibers of $p:\tYX\to
Y$. From the definition of the $V$-filtration it is clear that for any
$P\in gr^k_V\DX$ we have $\gvt P=P(\gvt - k)$; hence if $\CM$ is a
coherent $\DX$-module we may define an
endomorphism $\gTh$  of $gr_V\CM$ commuting with the action of $p_*\DtYXa$
by $\gTh=\gvt + k$ on $gr^k_V\CM$.

\begin{dfn}
A coherent $\DX$-module is said to be \textsl{specializable}
along $Y$ if, locally on $Y$, there exists a polynomial $b$ such that $b(\gTh)$
annihilates $gr_V\CM$.

The set of polynomials $b$ annihilating $gr_V\CM$ on an open set $U$
of $Y$ is an ideal of
$\C[T]$ and the generator of this ideal is called the
$b$\textsl{-function} for $\CM$ along $Y$ on $U$.
\end{dfn}

This $b$-function depends on the choice of the $V$-filtration but
its roots are independent of the choice of the $V$-filtration on
$\CM$ up to translations by integers.  If $\CI$ is a coherent
ideal of $\DX$, the $\DX$-module $\CM=\DX/\CI$ is provided with
the canonical $V$-filtration induced by the $V$-filtration of
$\DX$ and \textsl{then} the $b$-function of $\CM=\DX/\CI$ is
canonically defined.

In the same way, if $\CM$ is specializable and $u$ is a section of
$\CM$, the submodule $\DX u$ of $\CM$ is specializable \cite{ENS} and
it has a canonical $V$-filtration given by $(V^k\DX u)_{k\in \Z}$,
hence the $b$-function of $u$ is canonically defined.

Let $\gth$ be any differential operator on $X$ whose class in
$gr^0_V\DX$ is $\gvt$. Then, by definition of the $b$-function, there
is an operator $P\in V_{-1}\DX$ such that
\begin{equation} \label{b-equ}
(b(\gth)+P)u=0
\end{equation}

Assume now that $\CM$ is a holonomic $\DX$-module. Then $p^{-1}gr_V\CM$ is a holonomic
$\DtYXa$-module and the sheaf $\End_\DtYXa(p^{-1}gr_V\CM)$ is  (locally)
finite dimensional. Thus the endomorphism $\gTh$ has (locally) a
minimal polynomial which is, by definition, a $b$-function of
$\CM$. This means that {\sl holonomic $\DX$-modules are
specializable along any submanifold $Y$}.

If $Y$ is a hypersurface with local coordinates such that
$Y=\ensemble{(x,t)\in X}{t=0}$, the equation \ref{b-equ} is written as~:
\begin{equation} \label{b-equ2}
\left(b(tD_t)+tQ(x,t,D_x,tD_t)\right)u=0
\end{equation}

If $Y$ has codimension $d$ greater than $1$, this equation is~:
\begin{equation} \label{b-equ3}
\Big(b(\scal t{D_t})+\sum_{i=1}^dt_iP_i(x,t,D_x,[tD_t])\Big)u=0
\end{equation}
where $\scal t{D_t}=\sum t_iD_{t_i}$ and $[tD_t]$ is the collection of
all operators $(t_iD_{t_j})_{i,j=1\dots d}$.

\begin{dfn} \label{def:mon}
A section $u$ is said to be $1$-{\sl specializable} (or to have a
"{\sl regular $b$-function"}) if it satisfies an equation \ref{b-equ}
with an operator $P$ whose order is less or equal to the degree of the
polynomial $b$.

The $b$-function is \textsl{"monodromic"} if $u$ satisfies an equation
\ref{b-equ} with $P=0$.
\end{dfn}

\begin{rmk}
A holonomic $\DX$-module has always $b$-functions but in general,
it has no regular $b$-function (except if the module is regular
holonomic \cite{KKHOUCH}).

A monodromic $b$-function is less usual. It is coordinate dependent,
more precisely it depends on an identification of a neighborhood of
$Y$ in $X$ and a neighborhood of $Y$ in $\tYX$, e.g. a fiber bundle
structure of $X$ over $Y$.
\end{rmk}

\begin{rmk}
Let $f:X\to \C$ be a holomorphic function. The $b$-function of $f$
is usually defined as the generator of the ideal of polynomials
satisfying an equation $b(s)f^s(x)=P(s,x,D_x)f^{s+1}(x)$. This
$b$-function appears as a special case of the previous definition
if we consider the holonomic $\DX$-module $\DX\gd(t-f(x))$. Then
the equation $b(s)f^s(x)=P(s,x,D_x)f^{s+1}(x)$ is formally
equivalent to the equation
$b(-D_tt)\gd(t-f(x))=tP(-D_tt,x,D_x)\gd(t-f(x))$.
\end{rmk}

\subsection{Quasi-$b$-functions}\label{sec:gbfunct}

In this paper, we will use a new kind of $b$-functions.  In fact, we
want to replace the Euler vector field by a vector field $\sum
n_it_iD_{t_i}$. For example, let $\gf:X\to X$ be defined in a
coordinate system $(x,t)$ by $(x,t)\mapsto (x,s_1=t_1^{n_1},\dots,$
$s_p=t_p^{n_p})$ for some positive integers $(n_1,\dots,n_p)$,
$Y=\{s=0\}$, $\tY=\{t=0\}$ and $\CM$ a holonomic $\DX$-module. It is
known that the inverse image $\gf^*\CM$ is a holonomic
$\CaD_{\tX}$-module, then $\gf^*\CM$ will have a $b$-function along
$\tY$ and by direct image we will get a $b$-function for $\CM$ but
with $\sum s_iD_{s_i}$ replaced by $\sum n_is_iD_{s_i}$.

So, let us consider the fiber bundle $p:\tYX\to Y$. The sheaf $\DtYXr$
of relative differential operators is the subsheaf of
$\DtYXa$ of the differential operators on $\tYX$ which commute with all
functions of $p^{-1}\OY$. A differential operator $P$ on $\tYX$
is homogeneous of degree $0$ if for any function $f$ homogeneous
of degree $k$ in the fibers of $p$, $Pf$ is homogeneous of degree
$k$.

In particular, a vector field $\tgh$ on $\tYX$ which is a relative
differential operator homogeneous of degree $0$ defines a morphism
from the set of homogeneous functions of degree $1$ into itself which
commutes with the action of $p^{-1}\OY$, that is a section of
$$\HOM_{p^{-1}\OY}\left(\CO_\tYX[1],\CO_\tYX[1]\right)$$ and thus an
endomorphism of the dual fiber bundle $\TYX$.

Let $(x,t)$ be coordinates of $X$ such that $Y=\ensemble{(x,t)\in
X}{t=0}$. Let $(x,\tt)$ be the corresponding coordinates of
$\tYX$. Then $\tgh$ is written as~:
$$\tgh=\sum a_{ij}(x)\tt_iD_{\tt_j}$$ and the matrix $A=(a_{ij}(x))$
is the matrix of the associated endomorphism of $\CO_\tYX[1]$ which is
a locally free $p^{-1}\OY$-module of rank $d=codim_XY$. Its
conjugation class is thus independent of the choice of coordinates
$(x,t)$, as well as its eigenvalues which will be called the
eigenvalues of the vector field $\tgh$.

\begin{dfn}
A vector field $\tgh$ on $\tYX$ is {\sl definite
positive with respect to $Y$} on $U\subset Y$ if it is a relative
differential operator homogeneous of degree $0$ whose eigenvalues are
strictly positive rational numbers and which is locally diagonalizable as an
endomorphism of $\CO_\tYX[1]$.
We denote by $Tr(\tgh)$ the trace of $\tgh$.
\end{dfn}

A structure of \textsl{local fiber bundle of $X$ over $Y$} is  an
analytic isomorphism between a neighborhood of $Y$ in $X$ and a
neighborhood of $Y$ in $\tYX$. For example a local system of
coordinates defines such an isomorphism.

\begin{dfn}\label{def:defpos}
A vector field $\gh$ on $X$ is \textsl{definite positive with respect
to $Y$} if:

i)\ $\gh$ is of degree $0$ for the $V$-filtration associated to $Y$ and
the image $\gs_Y(\gh)$ of $\gh$ in $gr^0_V\DX$ is {\sl definite
positive with respect to $Y$} as a vector field on $\tYX$.

ii)\ There is a structure of local fiber bundle of $X$ over $Y$ which
identifies $\gh$ and $\gs_Y(\gh)$.

The eigenvalues and the trace of $\gh$ will be, by definition, those
of $\gs_Y(\gh)$.
\end{dfn}

It is proved in \cite[proposition 5.2.2]{ENS} that in the case where
$\gs_Y(\gh)$ is the Euler vector field $\gvt$ of $T_YX$ the condition (ii) is
always satisfied and that the local fiber bundle structure of $X$ over $Y$ is
essentially unique for a given $\gh$, but this is not true in general.

We will now assume that $X$ is provided with such a vector field
$\gh$. Let $\gb=a/b$ the rational number with minimum positive
integers $a$ and $b$ such that the eigenvalues of $\gb^{-1}\gh$ are
positive relatively prime integers. Let $\DX[k]$ be the sheaf of
differential operators $Q$ satisfying the equation $[Q,\gh]=\gb kQ$
and let $V^\gh_k\DX$ be the sheaf of differential operators $Q$ which
are equal to a finite sum (algebraic case) or a convergent series
(analytic case) $Q=\sum_{l\le k}Q_l$ with $Q_l$ in $\DX[l]$ for each
$l\in\Z$. This defines a filtration of $\DX$.

\begin{dfn}\label{def:quasi-b}
Let $u$ be a section of a $\DX$-module $\CM$. A polynomial $b$ is a
{\sl quasi-b-function} with respect to $\gh$ (or a $b(\gh)${\sl
-function} for short) if there exists a differential operator $Q$ in
$V^\gh_{-1}\DX$ such that $(b(\gh)+Q)u=0$.

The $b(\gh)$-function will be said to be  \textsl{regular} if the order of
$Q$ as a differential operator is less or equal to the order of
the polynomial $b$ and \textsl{monodromic} if $Q=0$.
\end{dfn}

If $\gs_Y(\gh)$ is the Euler vector field of $T_YX$, this definition
is essentially equivalent to the definition of the previous
section. However, the $V$-filtration is defined on the sheaf $\DX|_Y$
of differential operators defined in a neighborhood of $Y$ while, for
a given vector field $\gh$, the $V^\gh$-filtration is defined on any
open set where $\gh$ is defined. That is why it will be useful to
consider the second definition even in the case of the Euler vector
field.

If $\gh$ is given, we may locally diagonalize $\gs_Y(\gh)$ and
identify $\gh$ with $\gs_Y(\gh)$, that is assume that $\gh=\sum
n_it_iD_{t_i}$ and we may assume that the $n_i$ are integers after
multiplication of $\gh$ by an integer. In this case, the direct image
by the ramification $\gf$ associated to the $n_i$ of a $b$-function
for $\gf^*\CM$ is a $b(\gh)$-function for $\CM$, hence such a
$b(\gh)$-function always exists locally for holonomic $\DX$-module.

\subsection{Tame $\DX$-modules} \label{sec:tameres}
We will say that a cyclic holonomic $\DX$-module $\CM=\DX/\CI$\ \
is \textsl{tame} along a locally closed submanifold $Y$ of $X$ if
the roots of the $b$-function of $\CM$ relative to $Y$ are
strictly greater than the opposite of the codimension of $Y$. In
fact we will extend this definition by replacing the $b$-function
by quasi $b$-functions and also by introducing a parameter $\gd$.

\begin{dfn}\label{def:tamealong}
Let $\CM=\DX/\CI$ be a cyclic holonomic $\DX$-module and $Y$ be a
locally closed submanifold $Y$ of $X$. Let $\gd$ be a strictly
positive real number.

The module $\CM$ is $\gd$-\textsl{tame along $Y$} if $Y$ is open
in $X$ or if there exists a vector field $\gh$ on $X$ which is
definite positive with respect to $Y$ and a $b(\gh)$-function for
$\CM$ whose roots are all strictly greater than $-Tr(\gh)/\gd$.

The module $\CM$ is \textsl{tame} along $Y$ if it is $\gd$-tame
for $\gd=1$.
\end{dfn}

Let $\gh$ be a vector field on $X$ which is definite positive with
respect to a submanifold $Y$ (definition \ref{def:defpos}). A
subvariety of $X$ is conic for $\gh$ if it is invariant under the
flow of $\gh$, that is given by equations $(f_1,\dots,f_l)$
satisfying $\gh f_i=k_if_i$ for some integers $k_1,\dots,k_l$.

\begin{dfn}\label{def:conitame}
The cyclic module $\CM$ is \textsl{conic-tame} (resp.
$\gd$-\textsl{conic-tame)} \textsl{along $Y$} if $Y$ is open in
$X$ or if there exists a vector field $\gh$ on $X$ which is
definite positive with respect to $Y$ such that:

(i) there is a $b(\gh)$-function for $\CM$ whose roots are all
strictly greater than $-Tr(\gh)$ (resp. $-Tr(\gh)/\gd$)

(ii) the singular support of $\CM$ is conic for $\gh$.
\end{dfn}

Let us recall that the \textsl{singular support} of $\CM$ is set
of points of $X$ where its characteristic variety $\CHM$ is not
contained in the zero section of $T^*X$. If $\CM$ is holonomic,
its singular support is a nowhere dense subvariety of $X$.

\begin{rmk} \label{geom-monodrom} If $\CM$ admits a monodromic (quasi-)$b$-function,
the sections of $\CM$ are solutions of $b(\gh)u=0$ and the
characteristic variety of $\CM$ is contained in the subset of
$T^*X$ defined by $\gh=0$, this implies that the singular support
of $\CM$ is conic for $\gh$.
\end{rmk}

A stratification of the manifold $X$ is a union $X=\bigcup_\ga
X_\ga$ such that

\begin{itemize}
\item For each $\ga$, $\overline{X}_\ga$ is a complex algebraic (analytic) subset
  of $X$ and $X_\ga$ is its regular part.
\item $\{X_\ga\}_\ga$ is locally finite.
\item $X_\ga\cap X_\gb =\emptyset$ for $\ga\neq \gb$.
\item If $\overline{X}_\ga\cap X_\gb \neq \emptyset$ then $\overline{X}_\ga\supset X_\gb$.
\end{itemize}

If $\CM$ is a holonomic $\DX$-module, its characteristic variety
$\CHM$ is a homogeneous lagrangian subvariety of $\TX$ hence there
exists a stratification $X=\bigcup X_\ga$ such that
$\CHM\subset\bigcup_\ga T^*_{X_\ga}X$ \cite[Ch. 5]{KASHCOUR}.

\begin{dfn}\label{def:tame}
A cyclic holonomic $\DX$-module $\CM=\DX/\CI$, is \textsl{tame}
(resp. \textsl{conic-tame}, $\gd$-\textsl{tame},
\textsl{conic-$\gd$-tame}) if there is a stratification $X=\bigcup
X_\ga$ of $X$ such that $\CM$ is tame (resp. conic-tame,
$\gd$-tame, conic-$\gd$-tame) along any stratum $X_\ga$.
\end{dfn}

In the next sections, we will find $\DX$-modules which satisfy a
weaker condition:

\begin{dfn}\label{def:wtame}
The module $\CM$ is \textsl{weakly tame} if there is a
stratification $X=\bigcup X_\ga$ of $X$ such that for all $\ga$,
one of the two following conditions is satisfied:

(i)$\CM$ is tame along $X_\ga$

(ii)for each point $x$ of $X_\ga$, $\gp_\ga^{-1}(x)\cap
T^*_{X_\ga}X$ is not contained in the characteristic variety of
$\CM$.

Here, $\gp_\ga$ is the projection $T^*X\to X$ and $T^*_{X_\ga}X$
is the conormal bundle to $X_\ga$.
\end{dfn}

Any tame $\DX$-module is clearly weakly tame. The following properties will be proved in
the section \ref{sec:tame}:

\begin{thm}\label{thm:support}
If the $\DX$-module $\CM$ is weakly tame then it has no quotients with
support in a hypersurface of $X$.
\end{thm}

\begin{thm}\label{thm:suppdistr}
Let $M$ be a real analytic manifold and $X$ be a complexification
of $M$. If the $\DX$-module $\CM$ is weakly tame then it has no
distribution solution with support in a hypersurface of $M$.
\end{thm}

As pointed in section \ref{sec:generators}, if $\CM$ is a cyclic
$\DX$-module $\DX/\CI$, the solutions are defined as the common
solutions of all operators in $\CI$. Here "solutions" means solutions
in the usual meaning and do not concern the $\Ext^k$ of the module
for $k>0$. Any distribution solution of $\CM$ is analytic  where $\CM$
is elliptic \cite{SKK} that is outside the singular support of $\CM$. So, if $\CM$ is
tame the distributions solutions of $\CM$ are uniquely characterized
by their restriction to the complementary of the singular support of $\CM$.

\vskip 0.4cm

\begin{thm}\label{thm:ldloc}
Let $\CM$ be a holonomic {\textrm conic-tame} $\DX$-module with singular
support $Z$. Then any multivalued holomorphic function on $X-Z$ with polynomial
growth along $Z$ which is solution of $\CM$ extends uniquely as a
$L^2_{loc}$ solution of $\CM$ on $X$.
\end{thm}

\begin{thm}\label{thm:luloc}
Let $M$ be a real analytic manifold and $X$ be a complexification of
$M$. Let $\CM$ be a holonomic {\textrm conic-tame} $\DX$-module whose
singular support $Z$ is the complexification of a real sub-variety of
$\CM$.

Then any distribution solution of $\CM$ on an open subset of $M$ is a
$L^1_{loc}$-function.
\end{thm}

These theorems as well as the following proposition will be proved
in section \ref{sec:tame}.

\begin{prop}\label{prop:deltaloc}
Let $\gd>0$ and $\CM$ a holonomic {\textrm conic-$\gd$-tame}
$\DX$-module with singular support $Z$. Then any multivalued
holomorphic function on $X-Z$ with polynomial growth along $Z$ which
is solution of $\CM$ extends uniquely as a $L^{2\gd}_{loc}$ solution
of $\CM$ on $X$.

If $X$ is the complexification of a real analytic manifold $M$ and
$Z$ is the complexification of a real sub-variety of $M$, then any
distribution on an open subset of $M$ solution of $\CM$ is the sum
of a singular part supported by $Z$ and of a
$L^\gd_{loc}$-function. If $\gd>1$, the singular part is $0$.
\end{prop}

\begin{rmk}
If $\CM$ is a regular holonomic $\CaD$-module as defined by
Kashiwara-Kawa\"\i \cite{KKREG},  all multivalued holomorphic
solutions of $\CM$ on $X-Z$ have polynomial growth along $Z$  \cite{KKREG} and all
hyperfunction solutions are distributions \cite{KASHREG}. In this case, theorems
\ref{thm:ldloc}, \ref{thm:luloc} and proposition \ref{prop:deltaloc}
apply to all multivalued holomorphic solutions and hyperfunction
solutions.
\end{rmk}

\subsection{The Harish-Chandra theorem}\label{sec:HC}

Let $G_\R$ be a real semisimple Lie group and $\kg_\R$ be its Lie
algebra.

\noindent An \textsl{invariant eigendistribution} $T$ on $G_\R$ is a distribution which
satisfies:
\begin{itemize}
\item $T$ is invariant under conjugation by elements of $G_\R$.
\item $T$ is an eigendistribution of every biinvariant
differential operator $P$ on $G_\R$, i.e. there is a scalar $\gl$
such that $PT=\gl T$.
\end{itemize}

The main example of such distributions is the character of
an irreducible representation of $G_\R$. A famous theorem of
Harish-Chandra asserts that all invariant eigendistributions are
$L^1_{loc}$-functions on $G_\R$ \cite{HC1}.

Let us now explain what is the Hotta-Kashiwara $\CaD$-module \cite{HOTTA}.
Let $G$ be a connected complex semisimple group with Lie algebra
$\kg$. The group acts on $\kg$ by the adjoint action hence on the
space $\CO(\kg)$ of polynomial functions on $\kg$ which is
identified to the symmetric algebra $S(\kg^*)$ of the dual space
$\kg^*$. By Chevalley theorem, the space $\CO(\kg)^G\simeq
S(\kg^*)^G$ of invariant polynomials is equal to a polynomial
algebra $\C[P_1,\dots,P_l]$ where $P_1,\dots,P_l$ are
algebraically independent polynomials of $\CO(\kg)^G$ and $l$ is
the rank of $\kg$. In the same way, the space $\CO(\kg^*)^G\simeq
S(\kg)^G$ of invariant polynomials on $\kg^*$ is equal to
$\C[Q_1,\dots,Q_l]$ where $Q_1,\dots,Q_l$ are algebraically
independent polynomials of $\CO(\kg^*)^G$.

The differential of the adjoint action on $\kg$ induces a Lie algebra
homomorphism $\gt:\kg\to \mathrm{Der}S(\kg^*)$ by:
$$(\gt(a)f)(x)=\frac d{dt}f\left(\exp(-ta).x\right)|_{t=0}
\quad\mathrm{for}\quad a\in\kg, f\in S(\kg^*), x\in\kg$$ i.e.
$\gt(a)$ is the vector field on $\kg$ whose value at $x\in\kg$ is $[x,a]$.

An element $a$ of $\kg$ defines also a vector field with constant
coefficients on $\kg$ by:
$$(a(D_x)f)(x)=\scal{a}{df}=\frac d{dt}f(x+ta)|_{t=0}\quad\mathrm{for}\quad f\in S(\kg^*), x\in\kg$$

By multiplication, this extends to an injective morphism from
$S(\kg)$ to the algebra of differential operators with constant
coefficients on $\kg$. We will identify $S(\kg)$ with its image
and denote by $P(D_x)$ the image of $P\in S(\kg)$.

For $\gl\in\kg^*$, the Hotta-Kashiwara module $\MFl$ is the
quotient of  $\Dg$ by the ideal generated by $\gt(\kg)$ and by the
operators $Q(D_x)-Q(\gl)$ for $Q\in S(\kg)^G$.

A result of Harish-Chandra \cite[lemma 24]{HC3} shows that there is an
equivalence between the invariant eigendistributions and the solutions
of the Hotta-Kashiwara modules. More precisely, there is an analytic
function $\gf$ on $\kg$, invertible on a neighborhood of $0$, such
that $u$ is an invariant eigendistribution if and only if $X\mapsto
\gf(X)u\left(exp(X)\right)$ is a distribution solution of a
Hotta-Kashiwara module.

The definition of Hotta-Kashiwara extends to:

\begin{dfn}\label{def:KH}
Let $F$ be a finite codimensional ideal of $S(\kg)^G$. We denote
by $\CI_F$ the ideal of $\Dg$ generated by $\gt(\kg)$ and by $F$,
this defines the coherent $\Dg$-module $\MF=\Dg/\CI_F$.
\end{dfn}

The filtration induced on $F$ by the filtration of the
differential operators is the same than the filtration of the
symmetric algebra $S(\kg)$. If $F$ is finite codimensional, its
graduate is a power of $S_+(\kg)^G$, the set of non constant elements
of $S(\kg)^G$, hence defines the nilpotent
cone $\kN(\kg^*)$ of $\kg^*$. The cotangent bundle $T^*\kg$ is
identified with $\kg\x\kg^*$ and $\kg^*$ to $\kg$ by the Killing
form, then if $\kN(\kg)$ is the nilpotent cone  of $\kg$, the
characteristic variety of $\MF$ is \cite[4.8.3.]{HOTTA}:
$$\ensemble{(x,y)\in\kg\x\kg}{[x,y]=0, y\in \kN(\kg)}$$
and $\MF$ is a holonomic $\Dg$-module.

In particular, if $\kgrs$ is the set of regular semisimple
elements in $\kg$, $\MF$ is elliptic on $\kgrs$ and the singular
support of $\MF$ is the algebraic variety
$\kg'=\kg\setminus\kgrs=\ensemble{x\in\kg}{\gD(x)=0}$ where $\gD$ is
defined as follows. If $n=\dim \kg$, we set for $x\in\kg$:
$$\det(t.1 - \ad x) = \sum_{i=0}^n(-1)^{n-i}p_i(x)t^i$$
where $\ad x$ is the adjoint action of $x$ on $\kg$ that is $\ad
x(z)=[x,z]$. The rank $l$ of $\kg$ is the smallest integer $r$
such that $p_{\,r}\not\equiv 0$ and $\gD(x)\equiv p_{\,l}(x)$ is
the equation of $\kg'$.

If $\kg_\R$ is a real semisimple Lie algebra and $\kg$ its
complexification, $\gD(x)$ is a polynomial on $\kg$ which is real
on $\kg_\R$ hence $\kg'$ is the complexification of $\kg'_\R$. The
invariant eigenfunctions on $\kg_\R$ are distributions
solutions of $\MF$, in particular they are analytic on
$\kg_\R\setminus \kg'_\R$. The main result of this section is~:

\begin{thm}\label{thm:main}
The holonomic $\Dg$-module $\MF$ is conic-tame.
\end{thm}

As a consequence  we get the theorem of Harish-Chandra:

\begin{cor}\label{cor:main}
\
  \begin{enumerate}
  \item There is no invariant eigendistribution supported by $\kg'_\R=\kg_\R\setminus\kgrs$.
  \item The invariant eigendistribution are $L^1_{loc}$ functions on
    $\kg_\R$ and analytic on $\kg_\R\cap\kgrs$.
  \item The module $\MF$ has no quotient supported by $\kg'=\kg\setminus\kgrs$.
  \end{enumerate}
\end{cor}

It has been proved in \cite{LBY} that the module $\MF$ is a regular
holonomic $\Dg$-module. This implies that the results of
corollary \ref{cor:main} are true for hyperfunction solutions as well
as for distributions.

For each $s\in\kg$ semi-simple, $\kg^s=\ensemble{x\in\kg}{[x,s]=0}$ is
a reductive Lie algebra. Let $d_s$ and $r_s$ be the dimension and the
rank of the semi-simple Lie algebra $[\kg^s,\kg^s]$. Let $u$ be
the minimum of $r_s/d_s$ over all semi-simple elements of $\kg$ and
$\gd(\kg)=\frac{1+u}{1-u}$.

For example, if $\kg=\ks\kl_n(\C)$, $\gd(\kg)=1+\frac2n$.

\begin{prop}\label{prop:prec}
The holonomic $\Dg$-module $\MF$ is $\gd$-tame for any $\gd<\gd(\kg)$.
\end{prop}
This implies that the distribution solutions of $\MF$ are
$L^{\gd}_{loc}$-functions for $\gd<\gd(\kg)$.

\subsection{Symmetric pairs}\label{sec:Sym}

Let $G$ be a connected complex reductive algebraic Lie group with Lie
algebra $\kg$. Fix a non-degenerate, $G$-invariant symmetric bilinear
form $\gk$ on the reductive Lie algebra $\kg$ such that $\gk$ is the
Killing form on the semisimple Lie algebra $[\kg,\kg]$. Fix an
involutive automorphism $\gs$ of $\kg$ preserving $\gk$ and set
$\kk=\Ker(\gs-I)$, $\kp=\Ker(\gs+I)$.  Then $\kg=\kk\oplus\kp$ and
the pair $(\kg,\kk)$ or $(\kg,\gs)$ is called a symmetric pair. As
$\kk$ is a reductive Lie subalgebra of $\kg$, it is the Lie algebra of a
connected reductive subgroup $K$ of $G$. This group $K$ acts on $\kp$
via the adjoint action and the differential of this action induces a
Lie algebra homomorphism $\gt:\kk\to\textrm{Der}S(\kp^*)$ by:
$$(\gt(a).f)(x)=\frac d{dt}f\left(exp(-ta).x\right)|_{t=0} \textrm{ for
  }a\in\kk,f\in S(\kp^*),x\in \kp$$
(see \cite{LEVASS2} for the details)

If the group is $G\x G$ for some semisimple group $G$ and $\gs$
is given by $\gs(x,y)=(y,x)$, then $K\simeq G$, $\kp\simeq\kg$
and we find the definitions of section \ref{sec:HC}. We will call this
case the ``diagonal case''. In fact, the previous section is a special
case of this one but we give the definitions in the two cases for the
reader's convenience.

Let $x=s+n$ be the Jordan decomposition of $x\in\kp$, that is $s$
is semi-simple, $n$ is nilpotent and $[s,n]=0$. As this
decomposition is unique, if $x\in\kp$ then $s$ and $n$ are both
elements of $\kp$. The element $x$ is said to be \textsl{regular}
if the codimension of its orbit $K.x$ is minimal, this minimum is
the rank of the pair $(\kg,\kk)$ or of $\kp$. The set
$\kprs=\kgrs\cap\kp$ of semisimple regular elements of $\kp$ is
Zariski dense and its complementary $\kp'$ is defined by a
$K$-invariant polynomial equation $\gD(x)=0$ \cite{KOSR}. If
$(\kg,\kk)$ is a complexification of a real symmetric pair, this
equation is real on the real space.

The set $\kN(\kp)$ of nilpotent elements of $\kp$ is a cone.  If
$\CO_+(\kp)^K$ is the space of non constant invariant polynomials
on $\kp$, then $\kN(\kp)$ is equal to:
$$\kN(\kp)=\ensemble{x\in\kp}{\forall P\in\CO_+(\kp)^K\ \  P(x)=0}.$$
The set of nilpotent orbits is finite and define a stratification
of $\kN(\kp)$ \cite{KOSR}. By an extension of the Chevalley
theorem, the space $\CO(\kp)^K$ is a polynomial algebra
$\C[P_1,\dots,P_l]$ where $P_1,\dots,P_l$ are algebraically
independent polynomials of $\CO(\kp)^K$ and $l$ is the rank of
$\kp$.

If $F$ is a finite codimensional ideal of $\CO(\kp^*)^K=S(\kp)^K$, the module
$\MF$ is the quotient of $\Dp$ by the ideal generated by $\gt(\kk)$ and by
$F$. $\MF$ is a holonomic $\Dp$-module \cite{LEVASS2} whose characteristic
variety is contained in:
$$\ensemble{(x,y)\in\kp\x\kp}{[x,y]=0, y\in \kN(\kp)}$$
Its singular support is $\kp'=\kp\setminus\kprs$.

We proved in \cite{LBY} that the module $\MF$ is always regular (hence
the hyperfunctions solutions are distributions) but in some cases it
has non zero solutions with support a hypersurface (see \cite{SEKI} or
\cite{LEVASS2} for an example) hence is not always tame. We will show
that the roots of the $b$-functions are bounded below and the bounds will
be calculated from the numbers $\lpx$ defined by Sekiguchi \cite{SEKI}
\cite{LEVASS2}.

Let $x\in\kN(\kp)$, an extension of the Jacobson-Morosov theorem
\cite{KOSR} shows that there exists a
\textsl{normal S-triple} containing $x$, that is there exist $y\in\kp$ and
$h\in\kk$ such that $[x,y]=h$, $[h,x]=2x$ and $[h,y]=-2y$. Set $\ks=\C
h\oplus\C x\oplus \C y\simeq\sld(\C)$. The $\ks$-module $\kg$ decomposes as
$\kg=\bigoplus_{j=1}^sE(\gl_j)$, where $E(\gl_j)$ is a simple $\ks$-module of
highest weight $\gl_j\in\N$. We can choose a basis $(v_1,\dots,v_m)$ of
$\kp^y=\ensemble{z\in\kp}{[z,y]=0}$ with $[h,v_i]=-\gl_iv_i$ and we have $\kp=[x,\kk]\oplus\kp^y$. If
$\kg$ is semisimple, we define \cite{LEVASS2}:
\begin{equation}
\lpx=\sum_{j=1}^m(\gl_j+2)-\dim\kp\label{equ:lpx}
\end{equation}

\begin{rmk}
The number $m$ is the dimension of the space $\kp^y$ and $s$ is the
dimension of $\kg^y$. By proposition 5 in \cite{KOSR} we have $s=2m-k$
where $k=2\dim\kp-\dim\kg$ is independent of $y$.
\end{rmk}

Recall that a non zero nilpotent $x\in\kp$ is said
$\kp$-\textsl{distinguished} if $\kp^x\subset\kN(\kp)$. As the number
of nilpotent orbits is finite, there is only a finite number of
distinct numbers $\lpx$ for $x\in\kN(\kp)$ and we set:
$$\lp=\min\ensemble{\lpx}{x\in\kN(\kp), x\ \mathrm{distinguished}}$$
If $\kg$ is reductive but not semisimple, we consider $\tkg=[\kg,\kg]$,
$\tkk=\kk\cap\tkg$ and $\tkp=\kp\cap\tkg$. Then $(\tkg,\tkk)$ is a symmetric
pair with $\tkg$ semisimple and we set $\lp=\gl_{\tkp}$, we will also
consider the ``reduced dimension'' of $\kp$ that is $\redim\kp=\dim\tkp$.

Let $s\in\kp$ be semisimple, then $\kg^s=\kk^s\oplus\kp^s$ is a
symmetric pair and we define:
$$\gm_\kp=\min\ensemble{\frac12(\gl_{\kp^s}-\redim\kp^s)}{s\in\kp, s\ \mathrm{semisimple}}$$
(this minimum is again taken over a finite set).

Let $x$ be a point of $\kp$ with Jordan decomposition $x=s+n$, we set
$\lpx=\gl_{\kp^s}(n)$

In section \ref{sec:stratificationbis}, we will define a finite
stratification of $\kp$ and to each stratum $\gS$ we will associate a
vector field $\gh_\gS$ definite positive with respect to $\gS$
such that $\kp\setminus\kprs$ is conic relatively to $\gh_\gS$. If
$x=s+n$ is the Jordan decomposition, the numbers $\gm_{\kp^s}$,
$\lpx=\gl_{\kp^s}(n)$ and $\redim\kp^s$ are independent of
$x\in\gS$. We will denote $\gm_\gS=\gm_{\kp^s}$ and
$t_\gS=(\gl_{\kp^s}(n)+\redim\kp^s)/2$.

The main result of this paper is the following which will be
proved in section \ref{sec:proof}:
\begin{thm}\label{thm:mainsym}
The space $\kp$ admits a finite stratification and to each stratum $\gS$ is
associated a vector field $\gh_\gS$ definite positive with respect to $\gS$
and such that $\kp\setminus\kprs$ is conic relatively to $\gh_\gS$.  The trace
of $\gh_\gS$ is $t_\gS$.

Let $F$ be a finite codimensional ideal of $S(\kp)^K$. The holonomic
$\Dp$-module $\MF$ admits a $b(\gh_\gS)$-function along each stratum $\gS$
whose roots are greater or equal to $\gm_\gS$.
\end{thm}

\begin{rmk}
In the proof, we will see that the polynomial $b$ depend only on the
semisimple part of any $x\in \gS$, in particular $b$ is the same for
all nilpotent orbits.
\end{rmk}

\begin{cor}\label{cor:wpremain}
If $(\kg,\kk)$ is a symmetric pair such that for any $s\in\kp$ semisimple,
$\gl_{\kp^s}>0$, then for any $F$ finite codimensional ideal of $S(\kp)^K$, the
holonomic $\Dp$-module $\MF$ is weakly tame.
\end{cor}

\begin{cor}\label{cor:premain}
If $(\kg,\kk)$ is a symmetric pair such that for any $x\in\kp$,
$\lpx>0$, then for any $F$ finite codimensional ideal of $S(\kp)^K$, the
holonomic $\Dp$-module $\MF$ is conic-tame.
\end{cor}

The difference between the two corollaries is that in the first one we
ask that $\lpx>0$ for  elements $x$ whose nilpotent part is
distinguished, while in the second
the condition $\lpx>0$ is required for all elements. In the diagonal
case, that is in the case of a semisimple Lie group $G$ acting on its
Lie algebra $\kg$, we have $\dim\kg=\sum_{i=1}^s(\gl_i+1)$ hence
if $x$ is nilpotent, $\lpx$ is equal to the codimension of the orbit of $x$. This
number is thus always positive and theorem \ref{thm:main} is a special
case of corollary \ref{cor:premain}.

In the general case, Sekiguchi defined in \cite{SEKI} a class of
symmetric pairs, called ``nice pairs'', for which he proved that
$\gl_{\kp^s}>0$ for any $s\in\kp$ semisimple. So, in the case of
nice pairs, the module $\MF$ is weakly tame.

Consider now a real symmetric pair $\kg_\R=\kk_\R\oplus\kp_\R$,
its complexification is a complex symmetric pair as defined before
(see \cite{SEKI} for the details). Then we may consider the
distributions or the hyperfunctions on $\kp_\R$ which are
solutions of the module $\MF$.

\begin{cor}\label{cor:mainsym}
Under the condition of corollary \ref{cor:wpremain} we have:
  \begin{enumerate}
  \item There is no solution of $\MF$ supported by $\kp'_\R=\kp_\R\setminus\kprs$.
  \item The module $\MF$ has no quotient supported by $\kp'=\kp\setminus\kprs$.

and under the condition of  corollary \ref{cor:premain} we have moreover:
  \item The distributions on $\kp_\R$ solutions of $\MF$ are $L^1_{loc}$ functions.
  \end{enumerate}
\end{cor}

The first point has been proved by Sekiguchi \cite{SEKI} and
Levasseur-Stafford \cite{LEVASS2}. The third point is new and may
be improved by:
\begin{cor}\label{cor:secsym}
Let $\gd(\kp)$ be the minimum of $-t_\gS/\gm_\gS$ over all strata
$\gS$. The module $\MF$ is $\gd$-tame for any $\gd<\gd(\kp)$.
\end{cor}

\begin{rmk}
If the roots of the $b$-functions of a module $\MF$ are not integers, as in
the example of Levasseur-Stafford \cite[Remarks after theorem 3.8]{LEVASS2},
the module will satisfy points 1) and 2) of corollary \ref{cor:mainsym} by
remark \ref{rem:notinteger}, but the solutions will not be $L^1_{loc}$ if the
module is not tame.
\end{rmk}

As in the previous section, the module $\MF$ is a regular holonomic
$\Dg$-module by \cite{LBY}, hence all these results are true for
hyperfunction solutions as well as for distributions.

\section{Tame $\CaD$-modules}\label{sec:tame}

\subsection{Polynomials and differentials}\label{sec:inver}

Consider $\C^d$ with coordinates $(t_1,\dots,t_d)$ and denote by
$D_{t_1},\dots,D_{t_d}$ the derivations
$\dsur{t_1},\dots,\dsur{t_d}$. Let $(n_1,\dots,n_d)$ be strictly
positive integers and $\gh=\sum_{i=1}^dn_it_iD_{t_i}$. If
$\ga=(\ga_1,\dots,\ga_d)$ is a multi-index of $\N^d$, we denote
$\abs\ga=\sum \ga_i$ and $\scal\ga n= \sum\ga_in_i$. For $N\ge 0$, let
$$A_N=\ensemble{a\in\N}{\exists \ga\in\N^d, \abs \ga=N, a=\scal n\ga}$$
and define a polynomial $b_N$ by
$$b_N(T)=\prod_{k=0}^{N-1}\prod_{a\in A_k}(T+\abs n+a)$$

\begin{lem}\label{lem:binom}
For any $N\ge 1$, the differential operator $b_N(\gh)$ is in the left ideal
of $\CaD_{\C^d}$ generated by the monomials $t^\ga$ for $\abs\ga=N$.
\end{lem}

\begin{proof}
We prove the lemma by induction on $N$.  If $N=1$, $b_N(\gh)=\gh+\abs
n=\sum n_it_iD_{t_i}+n_i=\sum n_iD_{t_i}t_i$ is in the ideal generated
by $(t_1,\dots,t_d)$.

Let us denote $b'_N(T)=\prod_{a\in A_N}(T+\abs n+a)$. We remark that
$t^\ga\gh=(\gh-\scal \ga n)t^\ga$ hence if $\abs\ga=N$, $t^\ga
b'_N(\gh)=b'_N(\gh-\scal \ga n)t^\ga=c(\gh)(\gh+\abs
n)t^\ga=c(\gh)(\sum n_iD_{t_i}t_i)t^\ga$ is in the left ideal
generated by $t^\ga$ for $\abs \ga=N+1$. As $b_{N+1}(T)=
b_N(T)b'_N(T)$, if the lemma is true for $N$ it is true for
$N+1$.
\end{proof}

\begin{rmk}The main property of $b_N$ is that its roots are all integers lower or
equal to $-Tr(\gh)$. If $\gh$ is the Euler vector field of $\C^d$ that
is $n_1=\dots=n_d=1$, then $b_N(T)=(T+d)(T+d+1)\dots(T+d+N-1)$.
\end{rmk}

\begin{cor}\label{cor:support}
Let $Y$ be a smooth subvariety of $X$ of codimension $d$, $\CM$ be a
coherent $\DX$-module and $u$ be a section of $\CM$ with support in
$Y$. For any vector field $\gh$ definite positive with respect to $Y$,
$u$ has a $b(\gh)$-function along $Y$ which roots are integers lower
or equal to $-Tr(\gh)$.
\end{cor}

\begin{proof}
Let $(x_1,\dots,x_{n-d},t_1,\dots,t_d)$ be local coordinates such
that $Y=\{(x,t)\in X|t=0\}$ and $\gh=\sum_{i=1}^dn_it_iD_{t_i}$.
If $u$ is supported by $Y$ then there exists some integer $M$ such
that $t_1^Mu=\dots=t_d^Mu=0$. Let $N=(M-1)d+1$, then for any
monomial $t^\ga$ such that $|\ga|=N$ we have $t^\ga u=0$ and the
result comes from lemma \ref{lem:binom}.
\end{proof}

\subsection{$\CaD$-modules supported by a submanifold}\label{sec:supp}

\begin{prop}
Let $Y$ be a smooth subvariety of $X$ of codimension $d$, $\CI$ be a coherent
ideal of $\DX$ and $\CM=\DX/\CI$. Assume that $\CM$ is specializable and that
all the integer roots of the $b$-function $b$ are strictly greater than $-d$,
then $\CM$ has no quotient with support in $Y$.
\end{prop}

\begin{proof}
Let $\CN$ be a quotient of $\CM$ supported by $Y$ and $u$ the image of $1$ in
$\CN$. Then $\CN$ is generated by $u$. But the $b$-function of $u$ has all
roots strictly greater than $-d$ or non integers by hypothesis and all roots
are integers less or equal to $-d$ from corollary \ref{cor:support}, hence this
$b$-function must be equal to $1$, thus $u=0$ and $\CN=0$.
\end{proof}

There is a similar result for quasi-$b$-functions which we will prove
now. Let us first recall that, if $Y$ is a submanifold of $X$ and
$i:Y\to X$, the sheaf $\DYX$ is the $(\DY,i^{-1}\DX)$-bimodule defined
as $\OY\ox_{i^{-1}\OX}i^{-1}\DX$ and the sheaf $\DXY$ is the
$(i^{-1}\DX,\DY)$-bimodule defined as
$i^{-1}\DX\ox_{i^{-1}\gO_X}\gO_Y$. Here $\gO_X$ is he sheaf of
differential forms of maximum degree on $X$. We will consider $\DYX$
and $\DXY$ as $\DX$-modules supported by $Y$.

\begin{lem} \label{lem:dyx}
Let $Y$ be a smooth subvariety of $X$, and $\gh$ be a vector field
on $X$ which is definite positive with respect to $Y$. Let
$b\in\C[T]$ be a polynomial, $Q$ be an operator in $V^\gh_{-1}\DX$
and $P=b(\gh)+Q$. Let $\CN$ be a left $\DY$-module and $\CN\,'$ be
a right $\DY$-module.

(i) If all the integer roots of $b$ are strictly greater than
$-Tr(\gh)$ then $P$ is an isomorphism of $\DXY\ox_\DY\CN$.

(ii) If all the integer roots of $b$ are strictly negative then
$P$ is an isomorphism of $\CN\,'\ox_\DY\DYX$.
\end{lem}

\begin{proof}
Let us fix local coordinates of $X$ such that $Y=\ensemble{(x,t)\in
  X}{t_1=\dots=t_d=0}$. The $\DX$-module $\DXY$ is the quotient of
$\DX$ by the left ideal generated by $t_1,\dots,t_d$ hence the sections
of $\DXY$ may be represented by finite sums:
$$u=\sum u_{\ga,\gb}(x)D_x^\gb\gd^{(\ga)}(t)$$
where $\gd^{(\ga)}(t)$ is the class of $D_t^\ga$ modulo  $t_1,\dots,
t_d$.

We may change the coordinates and assume that $\gh=\sum
n_it_iD_{t_i}$, we may also multiply $\gh$ by an integer and
assume that all $n_i$ are integer (this modifies the polynomial
$b$ but the condition that all the roots of $b$ are strictly
greater than $-Tr(\gh)$ remains). We have
$\gh\,\gd^{(\ga)}(t)=-(\scal n\ga+\abs n)\,\gd^{(\ga)}(t)$.

Let us assume first that $\CN=\DY$. The image of the
$V^\gh$-filtration of $\DX$ on $\DXY$ is the filtration by  $\scal
n\ga$ hence this filtration is only in positive degrees. So, to
prove that $P$ is bijective on $\DXY$ it is enough to show that
$b(\gh)$ is bijective on the graduate $gr_V\DXY$, that is on
homogeneous elements
$$u=\sum_{\scal n\ga =N}u_{\ga,\gb}(x)D_x^\gb\gd^{(\ga)}(t)$$
Decomposing $b$ into linear factors we have to show that $\gh+a$ is bijective
if $a<Tr(\gh)$ or if $a$ is not an integer which is clear.

We consider now a left $\DY$-module $\CN$ and define a filtration
by
$$V^\gh_k\left(\DXY\ox_\DY\CN\right)=\left(V^\gh_k\DXY\right)\ox_\DY\CN$$
As the $V^\gh$-filtration is trivial on $\DY$, this filtration is
compatible with the $\DX$-module structure. As $b(\gh)$ acts on
the graduate $gr_V\DXY\ox_\DY\CN$ by $b(\gh)(A\ox u)=b(\gh)(A)\ox
u$, this action is bijective and $P$ is an isomorphism of
$\DXY\ox_\DY\CN$.

Let us now consider the sheaf $\DYX$. It is the quotient of
$\DX$ by the \textsl{right} ideal generated by $t_1\dots t_d$ hence the sections
of $\DXY$ may still be represented by finite sums:
$$u=\sum u_{\ga,\gb}(x)D_x^\gb\gd^{(\ga)}(t)$$ where $\gd^{(\ga)}(t)$
is the class of $D_t^\ga$ modulo $t_1\dots t_d$ but now $\DX$ operates
on the right and we have
$\gd^{(\ga)}(t).t_iD_{t_i}=+\ga_i\gd^{(\ga)}(t)$ and the same
calculus shows that $\gh+a$ is bijective on $\DYX$ if $a>0$.
\end{proof}

\begin{prop}\label{prop:age}
Let $Y$ be a submanifold of $X$, and $\gh$ be a vector field on
$X$ which is definite positive with respect to $Y$. Let $\CI$ be a
coherent ideal of $\DX$ and $\CM=\DX/\CI$. Assume that $\CM$
admits a $b(\gh)$-function whose integer roots are strictly
greater than $-Tr(\gh)$, then $\CM$ has no quotient with support
in $Y$.
\end{prop}

\begin{proof}
Let $\CN$ be a quotient of $\CM$ supported by $Y$. Let $u$ be the
image in $\CN$ of the class of the operator $1$ in $\CM$. Then
$\CN$ is generated by $u$ which is annihilated by an operator
$P=b(\gh)+Q$ with $Q\in V_{-1}\DX$ and the integer roots of $b$
are strictly greater than $-Tr(\gh)$. On the other hand, if $\CN$
is supported by $Y$ there exists a coherent $\DY$-module $\CN_0$
such that $\CN$ is isomorphic to $\DXY\ox_\DY\CN_0$ as a
$\DX$-module \cite{INV2}. Applying lemma \ref{lem:dyx} we get that
$P$ is an isomorphism of $\DXY\ox_\DY\CN_0$.
\end{proof}

By definition, the inverse image by $i$ of a
$\DX$-module $\CM$ is
$$\CM_Y=\DYX\otimes^\L_{i^{-1}\DX}i^{-1}\CM$$

It is known \cite{LSCHAP} that if $\CM$ is a specializable, $\CM_Y$ is
a complex of $\DY$-modules with coherent cohomology. If no root
of the $b$-function is an integer (this does not depend of a
generator of $\CM$), then $\CM_Y=0$. In this case, $\CM$ has no
quotient and no submodule supported by $Y$.

Lemma \ref{lem:dyx} means that $\CM=\DX/\DX P$ satisfies $\CM_Y=0$ if
$P=b(\gh)+Q$ and all the roots of $b$ are strictly negative. With the same
proof than Proposition \ref{prop:age} we deduce:

\begin{prop}\label{prop:inverse}
Let $Y$ be a smooth subvariety of $X$, and $\gh$ be a vector field
on $X$ which is definite positive with respect to $Y$. Let $\CI$
be a coherent ideal of $\DX$ and $\CM=\DX/\CI$. Assume that $\CM$
admits a $b(\gh)$-function whose integer roots are strictly
negative then the first cohomology group of $\CM_Y$, that is
$\CM_Y^0=\DYX\otimes_{i^{-1}\DX}i^{-1}\CM$, is equal to $0$.
\end{prop}

We consider now a real analytic manifold $M$ and a complexification
$X$ of $M$. The sheaf of distributions on $M$ will be denoted by
$\DbM$.

\begin{prop}\label{prop:distrib}
Let $Y$ be a submanifold of $X$ and $\gh$ be a vector field on $X$
which is definite positive with respect to $Y$. Let $\CI$ be a
coherent ideal of $\DX$ and $\CM=\DX/\CI$. Assume that $\CM$ admits a
$b(\gh)$-function whose integer roots are strictly greater than
$-Tr(\gh)$.

Then $\CM$ has no distribution solution with support in $Y\cap M$.

Moreover, if $\CM$ admits a regular $b(\gh)$-function whose integer
roots are strictly greater than $-Tr(\gh)$, $\CM$ has no hyperfunction
solution with support in $Y\cap M$.
\end{prop}

\begin{proof}
Let $u$ be a distribution solution supported by $Y\cap M$. As
$Y\cap M$ is an analytic subset of $M$, we may assume that the
support of $u$ is contained in an analytic subset $N$ of $M$ and
prove the proposition by descending induction on the dimension of
$N$. So, we take a point $x$ of the regular part of $N$ and will
prove that $u$ vanishes in a neighborhood of $x$. We may thus
assume that $N$ is smooth and denote by $N_\C$ the
complexification of $N$ which is a complex submanifold of $Y$.

A distribution supported by $N$ is written in a unique way as
$$u=\sum_{|\ga|\le m} a_\ga(x)\gd^{(\ga)}(t)$$ where $a_\ga(x)$ is a
distribution on $N$ and $\gd^{(\ga)}(t)$ is a derivative of the
Dirac distribution   $\gd(t)$ on $N$. In fact, we have
$\gNDbM\simeq \CaD_{X\leftarrow N_\C}\otimes_{\CaD_{N_\C}}\DbN$ hence
$$\gNDbM\simeq \DXY\otimes_\DY\CaD_{Y\leftarrow N_\C}\otimes_{\CaD_{N_\C}}\DbN.$$

So, if an operator $P$ satisfies the conditions of the first part
of lemma \ref{lem:dyx}, it defines an isomorphism of $\gNDbM$.

In the present situation, there is a surjective morphism:
$$\DX/{\DX P}\arr \CM\arr 0$$
where $P$ satisfies these conditions, hence $\HOM_\DX(\CM,\gNDbM)=0$.

In the case of hyperfunctions, we have $\gG_N(\mathcal
B_M)=\DIXY\otimes_\DIY\CaD^\infty_{Y\leftarrow
N_\C}\otimes_{\CaD^\infty_{N_\C}} \mathcal B_N$ where $\DIY$ is
the sheaf of differential operators of infinite order on $Y$ and
$\DIXY=\CaD^\infty_X\otimes_\DX\DXY$. From \cite[theorem
3.2.1.]{LMO} applied to the dual of $\CM$ which is a Fuchsian
module as well as $\CM$ we know that:
$$\RHOM_\DX(\CM,\DIXY)\simeq\RHOM_\DX(\CM,\DXY)$$
and we conclude as before.
\end{proof}

If the characteristic variety does not contain the conormal to $Y$,
these results are true without condition on the $b$-function:

\begin{prop}\label{prop:pasdecar}
If there is no point $x\in Y$ such that the characteristic variety of
$\CM$ contains $\gp^{-1}(x)\cap T^*_YX$ then
$\CM$ has no quotient supported by $Y$ and no hyperfunction or
distribution solution on $M$ supported by $Y\cap M$.
\end{prop}

\begin{proof}These results are well known, let us briefly recall how
they can be proved.

Let $\CN$ be a quotient of $\CM$. Its characteristic variety $Ch(\CN)$
is contained in $\CHM$. On the other hand $Ch(\CN)$ is involutive
hence if it is contained in $\gp^{-1}(Y)$, $\gp^{-1}(x)\cap Ch(\CN)$
is void or contains $\gp^{-1}(x)\cap T^*_YX$ for any $x\in Y$. So
$Ch(\CN)$ is void and $\CN=0$.

Let $u$ be a hyperfunction supported by $Y\cap M$ and solution of
$\CM$. As in the proof of proposition \ref{prop:distrib}, we may
assume that $u$ is supported by a submanifold $N$ of $M$ whose
complexification $N_\C$ is contained in $Y$. Let $supp(u)\subset N$ be
the support of $u$, $SS(u)\subset T^*M$ be the singular spectrum of
$u$ as defined in \cite[chap. I]{SKK}and $\gp_\R:T^*M\to M$ the
projection. For each $x\in supp(u)$, $SS(u)$ contains
$T^*_NM\cap\gp_\R^{-1}(x)$ and by \cite[Cor 3.1.2 ch. III]{SKK}
$SS(u)\subset\CHM$ hence $u=0$
\end{proof}

\begin{prop}\label{prop:smoothprecise}
Let $Y$ be a smooth connected hypersurface of $X$ and $\gD$ be an equation of
$Y$. Let $\CM=\DX/\CI$ be a holonomic $\DX$-module whose singular
support is contained in $Y$. Assume that the roots of the $b$-function
of $\CM$ along $Y$ are $\ga_1,\dots,\ga_N$ with multiplicity
$n_1,\dots,n_N$.

Then there exists a neighborhood $W$ of $Y$ in $X$ such that the
holomorphic solutions of $\CM$ on $X\setminus Y$ with moderate growth
on $Y$ are of the form:
$$f(x)=\sum_{i=1}^N\sum_{j=0}^{n_i-1}f_{ij}(x)\gD(x)^{\ga_i}\log(\gD(x))^j$$
where the functions $f_{ij}$ are holomorphic functions on $W$.
Moreover, for each $i$, the function $f_{i,n_i-1}$ does not
vanish identically on $Y$ except if all  $f_{ij}$ with the same
$i$ are equal to $0$.

If the $b$-function is regular, we do not have to assume that $f$ has
moderate growth.

If all roots are strictly greater than $-1$ the solutions are
$L^2_{loc}$.
\end{prop}

\begin{rmk}
For a precise definition of multivalued $L^2_{loc}$ solution see
\cite[ch. IV]{BJORK2}.
\end{rmk}

\begin{proof}
We choose local coordinates $(y,t)$ of $X$ such that
$\gD\equiv t$ and we may choose an operator $P=b(tD_t)+tQ(y,t,D_y,tD_t)$ such that
$Pf=0$ and $b$ is the $b$-function of $\CM$.

As $\CM$ is elliptic on $X\setminus Y$ all solutions of $\CM$ are
holomorphic on $X\setminus Y$ and all solutions on an open subset of
$X\setminus Y$ extend uniquely to the whole of $X\setminus Y$ as
ramified functions around $Y$. As $\CM$ is holonomic the set of
solutions is locally finite dimensional on $X\setminus Y$ and thus the
solutions are of finite determination and may be written as finite
sums:
$$f(y,t)=\sum f_{ij}(y,t)t^{\gl_i}\log(t)^j$$
with $f_{ij}(y,t)$ holomorphic on  $X\setminus Y$ and
$\gl_{i}-\gl_{k}\notin\Z$ if $i\ne k$ (see \cite[ch. IV]{BJORK2}
for the details).

If $W$ is a neighborhood of $Y$ which can be identified with a
neighborhood of the zero section of $T_YX$, this decomposition is
unique on $W$ and the space of sums $t^{\gl}\sum_j
f_{j}(y,t)\log(t)^j$ for fixed $\gl$ is invariant by $\DX$, hence each
of these terms in $f$ is a solution on $W\setminus Y$. So,  we may now
assume that
$$f(y,t)=t^\gl\sum_{j=0}^{n-1} f_j(y,t)\log(t)^j.$$

If $f$ is a distribution, $f$ is in the Nilsson class, that is the
functions $f_j$ are meromorphic \cite{BJORK2}. In the same way,
if the $b$-function is regular, all solutions on $X\setminus Y$
extend uniquely as Nilsson class solutions \cite[theorem
3.2.11.]{INV}.

Thus we may adjust the number $\gl$ so that all $f_j(y,t)$ are
holomorphic and at least one of $f_j(y,0)$ is not identically $0$. Let
$m\le n$ the highest integer such that $f_{m-1}(y,0)\not\equiv0$. We write
$f_j(y,t)=\sum_{k\ge0} f_{jk}(y)t^k$. Then the equation $Pf=0$
gives:
$$b(tD_t)\sum_{j=0}^{m-1} f_{j0}(y)t^\gl\log(t)^j = 0$$
This implies that $\gl$ is a root of $b$ and that $m$ is less or equal
to the multiplicity of the root. If $m<n$ and $k_0$ is the valuation
of $f_{n-1}(y,t)$ at $t=0$, we still have
$$b(tD_t)f_{n-1,k_0}(y)t^{\gl+k_0}\log(t)^{n-1} = 0$$
and this implies that $\gl+k_0$ is another root of $b$ with multiplicity
at least $n$.

We have proved that the solutions are locally of the form
$$f(x)=\sum_{i=1}^N\sum_{j=0}^{n_i-1}f_{ij}(x)\gD(x)^{\ga_i}\log(\gD(x))^j$$
but the functions $f_{ij}$ are uniquely determined by $f$ and $\gD$,
as they are holomorphic and $Y$ is connected they cannot vanish
identically on some open subset of $Y$ hence
this formula is global on a neighborhood of $Y$.

If all roots are strictly greater than $-1$, it is clear that the
solutions are $L^2_{loc}$ (the hypersurface is complex hence of
real codimension $2$).
\end{proof}

\subsection{Application to Tame $\CaD$-modules}\label{sec:mainproofs}

We will now prove the results announced in section \ref{sec:tameres}.
\begin{proof}[Proof of theorems \ref{thm:support} and \ref{thm:suppdistr}]
Let $\CN$ be a non zero quotient of $\CM$ which is supported by a
hypersurface of $X$. If $Z$ is the singular support of $\CM$,
$\CM$ is locally isomorphic to some power $\CO_X^{\,l}$ on
$X\setminus Z$, hence $\CN$ is supported by $Z$.

Now we consider the stratification of definition \ref{def:wtame}.  Let
$d$ be the minimum of the codimension of the strata on which $\CN$ is
a non zero module. We can choose a point of a stratum of codimension
$d$ where $\CN$ is not $0$. But if we apply proposition \ref{prop:age}
or proposition \ref{prop:pasdecar} to this stratum we get a
contradiction.

Theorem \ref{thm:suppdistr} is proved exactly in the same way using
propositions \ref{prop:distrib} and \ref{prop:pasdecar}.
\end{proof}

\begin{rmk}\label{rem:notinteger}
This proof uses propositions \ref{prop:age} and \ref{prop:distrib} hence if
the roots of the $b$-functions are not integers, the result is still true even
if they are less than the codimension of the stratum.
\end{rmk}

\begin{proof}[Proof of theorem \ref{thm:ldloc}]
Let $f(x)$ be a multivalued holomorphic function on $X-Z$ solution of
$\CM=\DX/\CI$. The argument follows by induction on the codimension of
the strata of the stratification given in the hypothesis.

From proposition \ref{prop:smoothprecise}, we know that $f$ is
$L^2_{loc}$ on a neighborhood of the smooth part of $Z$.  Now, let $S$
be a stratum and $a$ a point of $S$. By definition of a
stratification, there is a neighborhood $V$ of $a$ where all strata
except $S$ are of codimension strictly lower than the codimension of
$S$, thus $f$ is $L^2_{loc}$ on $V-S$. We may assume that $V$ is
compact, then  $f$ is $L^2$ on $V-\gO$ for any neighborhood
$\gO$ of $S$. We want to prove that $f$ is $L^2$ on $V$ that is on $V-Z$
because $Z$ is negligible.

From the hypothesis, after shrinking $V$, there are
local coordinates on $V$ and integers $(n_1,\dots,n_d)$ such that:

a) $V\cap S=\ensemble{(x_1,\dots,x_p,t_1,\dots,t_d)\in V}{t=0}$

b) $Z$ is quasi-conic in $(t_1,\dots,t_d)$ with weights
$(n_1,\dots,n_d)$

c) $\CM$ admits a quasi-$b$-function with weights
$(n_1,\dots,n_d)$ and roots greater than $-\sum n_i$, that is
there exists a polynomial $b$ whose roots are $>-\sum n_i$ and a
differential operator $Q$ in $V^\gh_{-1}\DX$ such that
$b(\gh)+Q\in\CI$.

We decompose $V-Z$ into a finite number of simply connected
quasi-cones such that on each of them, there is one of the
coordinates $t_1,\dots,t_d$ which does not vanish. Let $\gG$ be
one of them and assume that $t_1$ is the non-vanishing coordinate.
Consider coordinates $(x,s)$ defined by $t_1=s_1^{n_1}$ and
$t_i=s_is_1^{n_i}$ for $i=2\dots d$. These coordinates are well
defined if we restrict $\arg t_1$ to $]-\pi,\pi]$ and $\arg s_1$
to $]-\pi/n_1,\pi/n_1]$. We have
$\gG=\ensemble{(x,s)}{(x,s_2,\dots,s_d)\in W, 0<\abs{s_1}\le \gd,
-\pi/n_1<\arg s_1\le\pi/n_1}$ for some set $W$.

In these coordinates, the vector field $\gh$ is equal to
$s_1D_{s_1}$ and $D_{t_i}=s_1^{-n_i}D_{s_i}$. The
$V^\gh$-filtration is now defined by $s_1D_{s_1}$, hence it is the
usual $V$-filtration relative to $\{s_1=0\}$. The operator $Q$ is
in $V^\gh_{-1}\DX$ hence the ideal $\CI$ contains an operator
\[b(s_1D_{s_1})+s_1Q(x,s,D_x,s_1D_{s_1},D_{s_2},\dots,D_{s_d})\]
that is of a $b$-function relative to the hypersurface $\{s_1=0\}$
whose roots are greater than $-\sum n_i$. So, we can apply
proposition \ref{prop:smoothprecise} and we find that $f$ has the
form:
$$f(x,s)=\sum_{i=1}^N\sum_{j=0}^{n_i-1}f_{ij}(x,s)s_1^{\ga_i}\log(s_1)^j$$
where the complex numbers $\ga_i$ are roots of $b$ and the
functions $f_{ij}(x,s)$ are holomorphic on $W\x\C$.

As $V-Z$ is the finite union of sets like $\gG$, it is enough to
show that $f$ is $L^2$ on $\gG$. The functions
$f_{ij}(x,s)s_1^{\ga_i}\log(s_1)^j$ are linear combinations of the
determinations of the multivalued function $f$ hence if $f$ is
$L^2_{loc}$ on $V-Z$, the same is true for each of them. So we
may assume that $f$ is equal to $f_0(x,s)s_1^{\ga}\log(s_1)^j$ with
$f_0(x,s)$ holomorphic on a set
$\ensemble{(x,s)}{(x,s_2,\dots,s_d)\in W, \abs{s_1}\le\gd}$ and
$\ga>-\sum n_i$. Now we have:
\begin{align*}
\norm{f}_{L^2(\gG)}^2 =&\int_{\gG}\abs{f(x,t)}^2dx\wedge\overline{dx}
\wedge dt\wedge\overline{dt} \\
=&n_1^2\int_{W\x\{\abs{s_1}<\gd\}}
\abs{f_0(x,s)}^2\abs{\log(s_1)}^{2j} \abs{s_1}^{2(\ga-1+\sum n_i)}
dx\wedge\overline{dx}\wedge ds\wedge\overline{ds}
\end{align*}
The function $f_0(x,s)$ is $L^2$ on $W\x\{\gep<\abs{s_1}<\gd\}$ for
any $\gep>0$ and holomorphic in $s_1$, $\ga+\sum n_i>0$ hence
$f_0(x,s)\log(s_1)^js_1^{\ga-1+\sum n_i}$ is $L^2$ on
$W\x\{\abs{s_1}<\gd\}$. So $\norm{f}_{L^2(\gG)}$ is finite and $f$ is
$L^2_{loc}$ at $a$ which ends the proof.
\end{proof}

\begin{proof}[Proof of theorem \ref{thm:luloc}]
By the hypothesis, the singular support $Z$
of $\CM$ is the complexification of the real variety $M\cap Z$
which is nowhere dense in $M$. Hence $\CM$ is elliptic on an open dense
subset of $M$ and the solution $u$ is analytic on this open set.

Therefore, $u$ extends to a ramified holomorphic solution $f$ of $\CM$
on $U-Z$ where $U$ is an open subset of $X$. In the following we
may replace $X$ by $U$ and assume that $X=U$. As $u$ is a
distribution, $f$ has moderate growth on a neighborhood of $M\cap Z$,
hence $f$ is a $L^2_{loc}$
function according to theorem \ref{thm:ldloc}.

If the given stratification of $Z$ is the complexification of a
real stratification, we may use the same proof than
\ref{thm:ldloc}. In the general case, theorem \ref{thm:luloc} is
the direct consequence of theorem \ref{thm:ldloc} and of the
following lemma.
\end{proof}

\begin{lem}
Let $M$ be a real analytic manifold, $X$ a complexification of $M$ and
$\gD$ a real analytic function on $M$ which extends to a holomorphic
function $\tgD$ on $X$. Let $L= \gD^{-1}(0)$ and $Z= \tgD^{-1}(0)$.

Let $f$ be a Nilsson class function on $X-Z$, then $f$ is $L^2_{loc}$
on $X$ if and only if the restriction of $f$ to $M$ extends to a
$L^1_{loc}$ function on $M$.
\end{lem}

\begin{proof}
Let us consider a resolution of singularities of $\gD$ which extends to
a complex resolution of $\tgD$, hence we have a real analytic manifold
$\tM$, a subvariety $\tL$ with normal crossing, a proper analytic map
$\gga:\tM\to M$ which is an isomorphism $\tM-\tL\to M-L$ and their
complexifications $\tX$, $\tZ$ and $\gga_\C$ with the same properties.

It is proved in \cite[Proposition 4.5.3.]{BJORK2} that $f$ is
$L^2_{loc}$ on $X$ if and only if $f_\circ\gga_\C$ is $L^2_{loc}$ on
$\tX$ and the same proof shows that $f|_M$ is $L^1_{loc}$ on $M$ if
and only if $(f|_M)_\circ\gga$ is $L^1_{loc}$ on $\tM$. So, we have to
prove the result on $\tM$. As the result is local, we may assume
that there are coordinates $(x_1,\dots,x_n)$ such that
$\tL=\{x_1\dots x_d=0\}$ and their complexification $(z_1,\dots,z_n)$
such that $\tZ=\{z_1\dots z_d=0\}$.

From \cite[Proposition 4.4.1.]{BJORK}, we know that the function
$f_\circ\gga_\C$ is equal to a finite sum $\sum
f_{\ga,k}(z)z^\ga(\log(z))^k$ with
$(\log(z))^k=(\log(z_1))^{k_1}\dots(\log(z_d))^{k_d}$ and
$z^\ga=z_1^{\ga_1}\dots z_d^{\ga_d}$.  We may assume that
$f_{\ga,k}(0)\ne 0$ and that the multi-indexes $(\ga,k)$ are all
different. Then $f_\circ\gga_\C$ is $L^2_{loc}$ on $\tX$ if and
only if $\Re \ga_i>-1$ for all $\ga$ appearing in the sum and all
$i=1,\dots,d$.  But $(f|_M)_\circ\gga$ is $L^1_{loc}$ on $\tM$
under the same condition, which proves the lemma.
\end{proof}

Proposition \ref{prop:deltaloc} is proved in the same way.

\section{Semisimple Lie Algebras and symmetric pairs}\label{sec:semi}

\subsection{Stratification of a semisimple algebra}\label{sec:stratification}

In this section, we will define the stratification which will be
used to prove that the Hotta-Kashiwara module is tame. This
stratification is well known, see \cite{ATI} for example.

Let $G$ be a connected complex semisimple algebraic Lie group with
Lie algebra $\kg$. The orbits in $\kg$ are the orbits of the adjoint
action of $G$.

An element $x$ of $\kg$ is said to be \textsl{semisimple} (resp.
\textsl{nilpotent}) if $\ad(x)$ is semisimple (resp $\ad(x)$ is
nilpotent). Any $x\in\kg$ may be decomposed in a unique way as
$x=s+n$ where $s$ is semisimple, $n$ is nilpotent and $[s,n]=0$
(Jordan decomposition). $x$ is said to be \textsl{regular} if the
dimension of its centralizer $\kg^x=\ensemble{y\in\kg}{[x,y]=0}$
is minimal, that is equal to the rank of $\kg$.

As pointed in \S\ref{sec:HC}, the set $\kgrs$ of semisimple regular elements
of $\kg$ is Zariski dense and its complementary $\kg'$ is defined by a
$G$-invariant polynomial equation $\gD(x)=0$. If $\kg$ is a complexification of a
real algebra $\kg_\R$, this equation is real on $\kg_\R$.

The set $\kN(\kg)$ of nilpotent elements of $\kg$ is a cone.  Let
$\CO_+(\kg)^G$ be the space of non constant invariant polynomials on
$\kg$, then $\kN(\kg)$ is equal to:
$$\kN(\kg)=\ensemble{x\in\kg}{\forall P\in\CO_+(\kg)^G\ P(x)=0}.$$
The set of nilpotent orbits is finite and define a stratification of
$\kN$ \cite[Cor 3.7.]{KOS1}.

We fix a Cartan subalgebra $\kh$ of $\kg$ and denote by $W$ the Weyl
group $W(\kg,\kh)$. The Chevalley theorem shows  that $\CO(\kg)^G$ is
equal to $\C[P_1,\dots,P_l]$ where $(P_1,\dots,P_l)$ are
algebraically independent invariant polynomials and $l$ is the rank of
$\kg$, that the set of polynomials on $\kh$ invariant under $W$ is
$\CO(\kh)^W=\C[p_1,\dots,p_l]$ where $p_j$ is the restriction to $\kh$
of $P_j$ and that the restriction map $P\mapsto P|_\kh$ defines
an isomorphism of $\CO(\kg)^G$ onto $\CO(\kh)^W$ \cite[\S
4.9.]{VARADA}. The space $\kh/W$ is thus isomorphic to $\C^l$.

Let $\gF=\gF(\kg,\kh)$ be the root system associated to $\kh$. For
each $\ga\in\gF$ we denote by $\kg_\ga$ the root subspace
corresponding to $\ga$ and by $\kh_\ga$ the subset
$[\kg_\ga,\kg_{-\ga}]$ of $\kh$ (they are all $1$-dimensional).

Let $\CF$ be the set of the subsets $P$ of $\gF$ which are closed and
symmetric that is such that $(P+P)\cap\gF\subset P$ and $P=-P$. For
each $P\in\CF$ we define $\khp=\sum_{\ga\in P}\kh_\ga$,
$\kgp=\sum_{\ga\in P}\kg_\ga$, $\khpo=\ensemble{H\in\kh}{\ga(H)=0
\textrm{ if } \ga\in P}$ and $\khpp =\ensemble{H\in\kh}{\ga(H)=0
\textrm{ if } \ga\in P,\ga(H)\ne0 \textrm{ if } \ga\notin P}$.

The following results are well-known (see \cite[Ch VIII, \S3]{BOU}):

a) $\kqp=\khp+\kgp$ is a semisimple Lie subalgebra of $\kg$ stable
under $\ad\kh$ and $\khpo$ is an orthocomplement of $\khp$ for the
Killing form, $\khp$ is a Cartan subalgebra of $\kqp$. The Weyl group
$W_P$ of $(\kqp,\khp)$ is identified to the subgroup $W'$ of $W$ of
elements whose restriction to $\khpo$ is the identity \cite[theorem
4.15.17]{VARADA}.

b) $\kh+\kgp$ is a reductive Lie subalgebra of $\kg$ stable under
$\ad\kh$. For any $s\in\khpo$, $\kh+\kgp\subset\kg^s$ and
$\khpp=\ensemble{s\in\khpo}{\kg^s=\kh+\kgp}$.

c) Conversely, if $s\in\kh$, there exists a subset $P$ of $\gF$
which is closed and symmetric such that $\kg^s=\kh+\kgp$. $P$ is
unique up to a conjugation by $W$.

To each $P\in\CF$ and each nilpotent orbit $\kO$ of $\kqp$ we
associate a conic subset of $\kg$
\begin{equation}
S_{(P,\kO)}=\bigcup_{x\in\khpp}G.(x+\kO)\label{def:strat}
\end{equation}
where $G.(x+\kO)$ is the union of orbits of points $x+\kO$.

\begin{prop}\label{prop:strat}
The sets $S_{(P,\kO)}$ define a finite stratification $\gS_\kg$ of $\kg$.
\end{prop}

This proposition is a special case of proposition
\ref{prop:stratp} and we refer to it for the proof. Let us
describe some of the strata:
\begin{itemize}
\item If $P=\emptyset$, there is one associated stratum which is
the set $\kgrs$ of all regular semisimple points of $\kg$.

\item If $P=\{-\ga,\ga\}$, $\kqp$ is isomorphic to $\sld$ and if
$\kO$ is the non-zero orbit of $\kqp$, $S_{(P,\kO)}$ is exactly
the stratum of codimension $1$ which is the smooth part of $\kg'$
the hypersurface of equation $\gD(x)=0$ (cf. \S\ref{sec:HC}).

\item If$P=\gF$, $\kqp=\kg$ and the strata $S_{(P,\kO)}$ are the
nilpotent orbits of $\kg$.
\end{itemize}

\subsection{Stratification of a symmetric pair}\label{sec:stratificationbis}

Let $(\kg,\kk)$ be a symmetric pair with $\kg=\kk\oplus\kp$. A
Cartan subspace $\ka$ of $\kp$ is a maximal abelian subspace of
$\kp$ consisting of semisimple elements. Its dimension coincide
with the rank $l$ of $\kp$. If $W=N_K(\ka)/Z_K(\ka)$ is the
associated Weyl group, the Chevalley restriction theorem gives an
isomorphism $\CO(\kp)^K\simeq \CO(\ka)^W$ and $\CO(\kp)^K$ is
equal to $\C[P_1,\dots,P_l]$ where $(P_1,\dots,P_l)$ are
algebraically independent invariant polynomials. We denote by $V$
the vector space
$V=\mathrm{Spec}(\CO(\kp)^K)\simeq\mathrm{Spec}(\C[P_1,\dots,P_l])$,
by $\varpi:\kp\to V$ the projection and  by $\varpi_0:\ka\to V$
its restriction.

For $\ga\in\ka^*$, we set $\kg_\ga=\ensemble{x\in\kg}{\forall
  h\in\ka,\ [h,x]=\ga(h)x}$. The restricted root space
$\gF=\gF(\kg,\ka)$ is the set of $\ga\in\ka^*$ such that
$\kg_\ga\ne0$. The dimension of $\kg_\ga$ is not necessarily $1$ as in
the diagonal case but we have the following results \cite{RICH}:

a)The Cartan subspaces of $\kp$ are all conjugated by $K$ and any
semisimple element of $\kp$ belongs to one of them.

b)Let $\km=Z_\kk(\ka)$ be the centralizer of $\ka$ in $\kk$, then:
$$\kg=\km\oplus\ka\oplus\bigoplus_{\ga\in\gF}\kg_\ga$$

c)Let $\kg_{[\ga]}=\kg_\ga\oplus\kg_{-\ga}$, then
$\dim\kg_\ga=\dim\kg_{[\ga]}\cap\kp=\dim\kg_{[\ga]}\cap\kk$ and, if
$\gF^+$is the set of positive roots for some total order on $\ka^*$, we
have:
$$\kp=\ka\oplus\bigoplus_{\ga\in\gF^+}\kg_{[\ga]}\cap\kp$$

Let $P\in\gF$ be symmetric and closed. Define $\kg_P=\sum_{\ga\in
P}\kg_\ga$, $\kap=[\kg_P,\kg_P]\cap\ka$,
$\kp_P=\kg_P\cap\kp=\sum_{\ga\in\gF^+\cap P}(\kg_{[\ga]}\cap\kp)$,
$\kapo=\ensemble{h\in\ka}{\forall\ga\in P,\ \ga(h)=0}$ and
$\kapp=\ensemble{h\in\kapo}{\forall\ga\notin P,\ \ga(h)\ne0}$.

Let $s\in\kapp$, then $\kg^s=\km\oplus\ka\oplus\kg_P$ and
$\kp^s=\ka\oplus\kp_P$. Conversely, let $s\in\ka$ be semisimple
and define $P=\ensemble{\ga\in\gF}{\ga(s)=0}$. Then $P$ is closed
and symmetric and $\kg^s=\km\oplus\ka\oplus\kg_P$,
$\kp^s=\ka\oplus\kp_P$.

The decomposition $\kg^s=\kk^s\oplus\kp^s$ defines a symmetric
pair with the same rank and the group acting on $\kp^s$ is $K^s$.

To each $P\subset\gF$ closed and symmetric and to each nilpotent orbit
$\kO$ of $\ka\oplus\kp_P$, we associate a conic subset of $\kp$:
\begin{equation}
S_{(P,\kO)}=\bigcup_{x\in\kapp}K.(x+\kO)\label{def:stratsym}
\end{equation}

\begin{prop}\label{prop:stratp}
The sets $S_{(P,\kO)}$ define a finite stratification $\gS_\kp$ of
$\kp$, that is:

(i) Each stratum $S_{(P,\kO)}$ is a smooth locally closed
subvariety of $\kp$.

(ii) The number of strata is finite.

(iii)  The union of all strata is equal to $\kp$.

(iv) The strata are mutually disjointed.

(v) If $S_1$ and $S_2$ are two strata such that
$S_1\cap\overline{S_2}\neq \emptyset$ then
$S_1\subset\overline{S_2}$.
\end{prop}

\begin{rmk}\label{rem:prod}
If $\kg$ is a reductive Lie algebra which is not semisimple, let
$\kc$ be the center of $\kg$ and set $\tkg=[\kg,\kg]$,
$\tkp=\kp\cap\tkg$, $\tkk=\kk\cap\tkg$, $\tka=\ka\cap\tkp$ and
$\kc_\kp=\kp\cap\kc$. Then $(\tkg,\tkk)$ is a symmetric pair with
$\tkg$ semisimple and $\tkg=\tkk\oplus\tkp$. We have
$\kg=\kc\oplus\tkg$, $\kp=\kc_\kp\oplus\tkp$ and
$\ka=\kc_\kp\oplus\tka$. Moreover $\gF(\kg,\ka)=\gF(\tkg,\tka)$
and a set $S_{(P,\kO)}$ in $\kp$ is exactly the direct sum of
$\kc_\kp$ and of the corresponding set $S_{(P,\kO)}$ in $\tkp$.
The stratification of $\kp$ is thus the direct sum of $\kc_\kp$
and of the stratification of $\tkp$.
\end{rmk}

Let $s\in\kp$ be semisimple. The previous remark apply to $\kp^s$:

 Set $\tkg_s=[\kg^s,\kg^s]$, $\tkp_s=\kp^s\cap\tkg$,
$\tkk_s=\kk^s\cap\tkg$. Then $(\tkg_s,\tkk_s)$ is a symmetric pair
with $\tkg_s$ semisimple and $\tkg_s=\tkk_s\oplus\tkp_s$. If
$\kc_s$ is the center of $\kg^s$,
$\kp^s=(\kp^s\cap\kc_s)\oplus\tkp_s$
\cite{LEVASS2}.

Remark also that $\kN(\kp^s)=\kN(\kp)\cap\kp^s$ and $a\in\kp^s$ is
semisimple in $\kg^s$ if and only if $a$ is semisimple in $\kg$.

\begin{lem}\label{lem:semi}
(i)There is a neighborhood $\gO$ of $s$ and an open embedding
$\gro:\gO\to \gro(\gO)\subset Ks\times\kp^s$ such that the
intersection of a $K$-orbit in $\kp$ with $\gO$ is equal to
$\gro^{-1}((Ks\times U)\cap\gro(\gO))$ where $U$ is a $K^s$-orbit
in $\kp^s$.

(ii)If  $x\in\kp^s\cap\gO$, $K.x\cap\kp^s\cap\gO$ is equal to
$K^s.x\cap\gO$
\end{lem}

This lemma is an easy consequence of  prop. \ref{prop:semi} and
its proof is postponed to section \ref{sec:proof}.

\begin{proof}[Proof of proposition \ref{prop:stratp}]
(i) The set $\kapp$ is transversal to the $K$-orbits of $\kp$,
hence each set $S_{(P,\kO)}$ is a smooth locally closed subvariety
of $\kp$.

(ii) As the set of subsets $P$ and the set of nilpotent orbits of
each $\kp^s$ are both finite, there is also a finite number of
sets $S_{(P,\kO)}$.

(iii) If $x\in\kp$ and $x=s+n$ is the Jordan decomposition, the
semisimple part $s$ belongs to a Cartan subspace of $\kp$. As they
are all conjugated we may assume $s\in \ka$. As remarked befoer,
there is some $P$ such that $\kp^s = \ka+\kp_P$. Then $n\in\kp_P$
and if $\kO$ is its $K^s$-orbit, $x\in S_{(P,\kO)}$. So $\kp$ is
the union of all $S_{(P,\kO)}$.

(iv) Let $S_1=S_{(P_1,\kO_1)}$ and $S_2=S_{(P_2,\kO_2)}$. If they
are not disjointed, there exist $k_1\in K$, $k_2\in K$,
$x_1\in(\ka_{P_1}^\bot)'$, $x_2\in(\ka_{P_2}^\bot)'$,
$n_1\in\kO_1$ and $n_2\in\kO_2$ such that
$k_1.(x_1+n_1)=k_2.(x_2+n_2)$. Remark that for $i=1,2$, $x_i$ is
semisimple, $n_i$ is nilpotent and $[x_i,n_i]=0$. By unicity of
the Jordan decomposition this implies $k_1.x_1=k_2.x_2$ and
$k_1.n_1=k_2.n_2$ hence $k=k_2^{-1} k_1$ is in the normalizer of
$\ka$ in $K$. Then $P_1=k.P_2$ and $\kO_1 =k.\kO_2$. As
$S_{(P,\kO)}=S_{(k.P,k.\kO)}$, we can conclude that $S_1=S_2$.

(v) Let $S_1=S_{(P_1,\kO_1)}$ and $S_2=S_{(P_2,\kO_2)}$ be two
distinct strata such that $S_1\cap\overline{S_2}\neq \emptyset$.
Let $x\in S_1\cap\overline{S_2}$. By definition of $S_1$,
$x=k(s+n)$ with $k\in K$, $s\in(\ka_{P_1}^\bot)'$ and $n\in\kO_1$.
We may replace $x$ by $s+n$ and assume $s+n\in
S_1\cap\overline{S_2}$.

The projection $\varpi_0:\ka\to V$ is closed and the set $\kapdo$
is closed, hence $\varpi^{-1}\varpi_0(\kapdo)$ is a closed subset
of $\kp$. If $x_2=k_2(s_2+n_2)\in S_2$,
$\varpi(x_2)=\varpi(k_2s_2)=\varpi_0(s_2)$ \cite[Lemma 12]{KOSR}
hence $\varpi^{-1}\varpi_0(\kapdo)$ contains $S_2$ and thus its
closure $\overline{S_2}$. Therefore
$s\in\varpi^{-1}\varpi_0(\kapdo)\cap\ka=\varpi_0^{-1}\varpi_0(\kapdo)$.
By definition of $\varpi_0$, this means that there exists $\gs\in
W$ such that $s^\gs\in\kapdo$. We may replace $x$ by $x^\gs$,
$P_1$ by ${P_1}^\gs$ and assume that $s\in\kapdo$. So we have
$s\in(\ka_{P_1}^\bot)'$ and $s\in\kapdo$, hence by definition of
$(\ka_{P_1}^\bot)'$ we have $P_2\subset P_1$. This shows in
particular that $(\ka_{P_1}^\bot)'\subset\kapdo$ and
$\kp_{P_2}\subset\kp_{P_1}$.

Now we apply lemma \ref{lem:semi}. As $S_1$ and $S_2$ are both
union of $K$-orbits, they are locally (i.e. in $\gO$) the product
of $K.s$ by their intersection with $\kp^s$. So
$\overline{S_2}\cap\kp^s=\overline{S_2\cap\kp^s}$ and to prove
that $S_1\subset \overline{S_2}$ it suffices to prove that
$S_1\cap\kp^s\subset\overline{S_2\cap\kp^s}$.

Now $\kO_1$ is a subset of $\kp^s$ by definition while
$\kp_{P_2}\subset\kp_{P_1}$implies that $\kO_2$ is also a subset
of $\kp^s$. Thus lemma \ref{lem:semi}(ii) shows that for $i=1,2$
$S_i\cap\kp^s=K^s.((\ka_{P_i}^\bot)'\oplus\kO_i)$. So we may now
replace $\kp$ by $\kp^s$ that is assume that $P_1=\gF$.

By remark \ref{rem:prod}, we may assume that $\kp$ is reduced, and
in this case with $P_1=\gF$, $S_1$ is a nilpotent orbit of $\kp$.
As $\overline{S_2}$ is invariant under the action of $K$, if a
point of $S_1$ is in  $\overline{S_2}$  then
$S_1\subset\overline{S_2}$.

This proof is valid only in the neighborhood $\gO$ of $s$. If
$x=s+n$ is a point of $S_1\cap\overline S_2$ but not in $\gO$, as
the set of nilpotent points of $\tkp$ is a cone, there is some
$k_0\in K^s$ such that $k_0.(s+n)=s+k_0.n$ is in $\gO$ and this
point is in $S_1\cap\overline S_2$. So $k_0^{-1}\gO$ is a
neighborhood of $x$ such that $S_1\cap k_0^{-1}\gO\subset\overline
S_2$.

This shows that $S_1\cap \overline S_2$ is an open subset of
$S_1$, as it is also a closed subset and $S_1$ is connected, we
have $S_1\cap\overline S_2=\emptyset$ or $S_1\subset\overline
S_2$.
\end{proof}

\subsection{Fourier transform}\label{sec:fourier}
We recall here a few things about the Fourier transform of
$\CaD$-modules after \cite{HOTTA}. Let $V$ be a finite-dimensional
vector space over $V$ and $\CaD_{[V]}$ the sheaf of algebraic
differential operators on $V$. Then
$$\gG(V,\CaD_{[V]})=\C[V]\ox_\C\C[V^*]=S(V^*)\ox_\C S(V)$$ where
$\C[V]=S(V^*)$ is the ring of polynomials on $V$ and
$\C[V^*]=S(V)$ is identified to the ring of constant coefficient
operators on $V$. In this way $\gG(V,\CaD_{[V]})$ is the
$\C$-algebra generated by $V\oplus V^*$ with the relations
$[v,w]=[v^*,w^*]=0$ and $[v,v^*]=-\scal{v}{v^*}$ for $v,w$ in $V$
and $v^*,w^*$ in $V^*$.

The Fourier transform is the isomorphism $\gG(V,\CaD_{[V]})
\to\gG(V^*,\CaD_{[V^*]})$ generated by the map $V\oplus V^*\to
V^*\oplus V$, $(v,v^*)\mapsto(v^*,-v)$. We denote by $\widehat P$ the
image of $P\in\gG(V,\CaD_{[V]})$ under this isomorphism.

The category of coherent $\CaD_{[V]}$-modules is equivalent to that of
finitely generated $\gG(V,\CaD_{[V]})$-modules, hence this define the
Fourier transform as a functor from the category of coherent
$\CaD_{[V]}$-modules onto the category of coherent
$\CaD_{[V^*]}$-modules. If $\CI$ is an ideal of $\CaD_{[V]}$ generated
by operators $P_1,\dots,P_N$ of $\gG(V,\CaD_{[V]})$, the Fourier
transform $\widehat\CI$ is generated by $\widehat P_1,\dots,\widehat
P_N$ and the Fourier transform of $\CM=\CaD_{[V]}/\CI$ is
$\widehat\CM=\CaD_{[V^*]}/\widehat\CI$. It is known that $\widehat\CM$
is holonomic if and only if $\CM$ is holonomic and the Fourier
transform of
$\widehat\CM$ is $a^*\CM$ with $a:V\to V$ given by $a(v)=-v$.

If we choose linear coordinates $(x_1,\dots,x_n)$ of $V$ and dual
coordinates $(\gx_1,\dots,\gx_n)$ of $V^*$, the Fourier transform is
given by $x_i\mapsto D_{\gx_i}$ and $D_{x_i}\mapsto -\gx_i$. If
$\gth=\sum x_iD_{x_i}$ is the Euler vector field of $V$ and
$\gth^*=\sum \gx_iD_{\gx_i}$ is the Euler vector field of $V^*$, we
have $\widehat\gth=-\gth^*-n$.

More generally, if $u=\sum u_{ij}x_iD_{x_j}$, then $\widehat u=-\sum
u_{ij}D_{\gx_i}\gx_j=-\sum u_{ij}\gx_jD_{\gx_i}-\sum u_{ii}$. So, if
$v:V\to V$ is a linear morphism, it defines a section $V\to TV=V\x V$
by $x\mapsto (x,v(x))$, that is a vector field $u$ on $V$ and if the
trace of $v$ is $0$, $\widehat u$ is the vector field associated to
the transpose ${}^tu:V^*\to V^*$.

A $\CaD_{[V]}$-module $\CM$ is \textsl{homogeneous} or
\textsl{monodromic} if it admits a monodromic $b$-function at
$\{0\}$, that is if for any section $u$ of $\CM$, there exists a
polynomial $b$ such that $b(\gth)u=0$.

\begin{prop}\label{prop:fourier}
Let $\CM=\CaD_{[V]}/\CI$ be a monodromic $\CaD_{[V]}$-module,
$\widehat\CM=\CaD_{[V^*]}/\widehat\CI$ its Fourier transform , $u$
the canonical generator of $\CM$ and $\widehat u$ the canonical
generator of $\widehat \CM$.

1)$\widehat\CM$ is monodromic and if $b$ is the $b$-function of
$u$ then $b(-\gth^*-n)$ is the $b$-function of $\widehat u$.

2) Let $x_0\in V$ and $b$ a polynomial such that $x_0$ is not in the
support of $b(\gth)u$, then the characteristic variety of
$\CaD_{[V^*]}b(-\gth^*-n)\widehat u$ does not meet $V\x \{x_0\}$.
\end{prop}

\begin{proof} The characteristic variety of $\CM$ is a subset of
$T^*V\simeq V\x V^*$ and the characteristic variety of $\widehat\CM$
is a subset of $T^*V^*\simeq V^*\x V$. If $x_0$ is not in the support
of $b(\gth)u$, there exists a function $a(x)$ such that $a(x_0)\ne0$
and $a(x)b(\gth)u=0$. As $\CM$ is monodromic, its support is conic and
we may assume that $a$ is a homogeneous function of $x$. Then
$a(D_\gx)b(-\gth^*-n)\widehat u=0$ which shows the proposition.
\end{proof}

Assume now that we are given a symmetric pair $(\kg,\kk)$ and that
$V=\kp$. The bilinear form $\gk$ defines an isomorphism
$\kp\simeq\kp^*$ which exchanges the morphism $\ad a:\kp\to\kp$
and its adjoint for any $a\in\kk$. With this identification, the
Fourier transform is an isomorphism of $\gG(\kp,\CaD_{[\kp]})$
onto itself. If $a\in\kk$, then $\gt(a)$ is the vector field on
$\kp$ defined by the linear morphism $\ad a$ whose trace is null
hence by what we said, its Fourier transform is $-\gt(a)$. On the
other hand, the isomorphism $\kp\simeq\kp^*$ extends to a
$K$-isomorphism $S(\kp)\simeq S(\kp^*)$ hence defines an
isomorphism $\tilde\gk:S(\kp)^K\simeq S(\kp^*)^K$. We get:

\begin{prop}\label{prop:fou}
Let $F$ be a finite codimensional ideal of $S(\kp)^K$, then
$\tilde\gk(F)$ is a finite codimensional ideal of
$S(\kp^*)^K=\C[\kp]^K$.

Let $\CM_{[F]}$ be the quotient of $\CaD_{[\kp]}$ by the ideal
generated by $\gt(\kk)$ and $F$. Its Fourier transform is the quotient
of $\CaD_{[\kp]}$ by the ideal generated by $\gt(\kk)$ and
$\tilde\gk(F)\subset\C[\kp]^K$.
\end{prop}

Let us denote by $\widehat{\CM_{[F]}}$ the Fourier transform of
$\CM_{[F]}$. The $\Dp$-module $\MF$ is by definition
$\Dp\otimes_{\CaD_{[\kp]}}\CM_{[F]}$ and we will denote by $\MFf$ the
module $\Dp\otimes_{\CaD_{[\kp]}}\widehat{\CM_{[F]}}$. We will call
this module the Fourier transform of $\MF$.

\subsection{Proof of main theorems}\label{sec:proof}

Consider a symmetric pair $(\kg,\kk)$ with $\kg=\kk\oplus\kp$, $\kg$ semi-simple, and a
nilpotent point $x\in\kp$. Let $\kO$ be the orbit of $x$ under the
action of $K$. As in section \ref{sec:Sym}, we consider a normal
$\sld$-triple  $(h,x,y)$ which defines a basis $(y_1,\dots,y_r)$ of
$\kp^y$ such that $[h,y_i]=-\gl_iy_i$. The number $\gl_\kp(x)$ is by definition
$\sum_{i=1}^r(\gl_i+2)-\dim\kp$ and we have
$\kp=[x,\kk]\oplus\kp^y$.

\begin{lem}\label{lem:coord}
There are local coordinates $(z_1,\dots,z_{n-r},t_1,\dots,t_r)$ of
$\kp$ near $x$ such that:

(i) $\kO=\ensemble{(z,t)}{t=0}$

(ii) The vector fields $D_{z_1},\dots,D_{z_{n-r}}$ are in the ideal of
$\Dp$ generated by $\gt(\kk)$.

(iii) The vector field
$$\gh_\kO=\sum_{i=1}^r(\frac{\gl_i}2+1)t_iD_{t_i}$$
is definite positive with respect to $\kO$ and its trace is equal to
$(\gl_\kp(x)+\dim\kp)/2$.

(iv) The Euler vector field $\gth$ of $\kp$ is equal to
$$\gth=\gh_\kO+\frac12D_{z_{n-r}}$$
\end{lem}

\begin{proof}
Let $x\in\kO$ and $(h,x,y)$ a normal $\sld$-triple. Let $V$ be a
linear subspace of $\kk$ such that $\kk=V\oplus\kk^x$ and $h$ is in
$V$. Let $b_1,\dots,b_{n-r}$ be a basis of $V$ with $b_{n-r}=h$. The
map $F:\C^{n-r}\x\C^r\to\kp$ given by
$$F(z_1,\dots,z_{n-r},t_1,\dots,t_r)=
\exp(z_1b_1)\dots\exp(z_{n-r}b_{n-r}).(x+\sum t_iy_i)$$ is a local
isomorphism hence defines local coordinates of $\kp$. In these
coordinates, $\kO$ is $\ensemble{(z,t)}{t=0}$ and  the vector fields
$D_{z_1},\dots,D_{z_{n-r}}$ are in the ideal of $\Dp$ generated by
$\gt(\kk)$ \cite[lemma 3.7.]{SEKI}.

The numbers $\gl_i$ are nonnegative integers hence $\gh_\kO$
is definite positive with respect to $\kO$ and its trace is equal to
$(\gl_\kp(x)+\dim\kp)/2$ by definition of $\gl_\kp(x)$.

By definition, the Euler vector field $\gth$ of $\kp$ acts as
$\gth(f)(u)=\frac d{ds}f(e^{s}u)|_{s=0}$ in linear coordinates $u$ of
$\kp$, hence in coordinates $(z,t)$:
$$\gth(f)(z,t)=\frac d{ds}f(F^{-1}(e^{s}F(z,t)))|_{s=0}$$
We have $[h,x]=2x$ and $[h,y_i]=-\gl_iy_i$ hence
$$\exp(sh).(x+\sum t_iy_i)=e^{2s}x+\sum e^{-\gl_i s}t_iy_i$$
therefore, as $b_{n-r}=h$ this gives:
$$F(z_1,\dots,z_{n-r-1},z_{n-r}+s/2,e^{s(\gl_1/2+1)}t_1,\dots,e^{s(\gl_r/2+1)}t_r)=
e^sF(z,t)$$
and thus
$$\gth=\sum_{i=1}^r(\frac{\gl_i}2+1)t_iD_{t_i}+\frac12D_{z_{n-r}}$$
\end{proof}

Let $F$ be a finite codimensional ideal of $S(\kp)^K$ and assume that
$F$ is graduate.  Let $\MFf$ be the Fourier transform of the module
$\MF$. By proposition \ref{prop:fou}, it is the quotient of $\Dp$ by
the ideal generated by $\gt(\kk)$ and by $\tilde\gk(F)$. As
$\tilde\gk(F)$ is graduate and of finite codimension it contains a
power of $\CO_+[\kp]^K=S_+(\kp^*)^K$ and $\MFf$ is supported by
$\kN(\kp)$.

Define $\tilde\gl_\kp$ as the minimum of $\gl_\kp(x)$ over all
nilpotents $x\in\kp$.

\begin{prop}\label{prop:grad}
The $b$-function of $\MFf$ at $\{0\}$ is monodromic and its roots
are lower or equal to $-(\tilde\gl_\kp+\dim\kp)/2$.
\end{prop}

\begin{proof}
The module $\MFf$ is supported by $\kN(\kp)$ which is a finite union
of nilpotent orbits \cite{KOSR}. By descending induction on the
dimension of these orbits, we have to prove that if $v$ is a section
of $\MFf$ supported by a nilpotent orbit $\kO$ in a neighborhood of a
point $x$ of this orbit, then there is a polynomial $b_x$ with roots
lower or equal to $-(\tilde\gl_\kp+\dim\kp)/2$ such that $b_x(\gth).v$
vanishes on a neighborhood of $x$ hence on a neighborhood of the orbit
of $x$.

So let $v$ be a section of $\MFf$ supported by $\kO$ in a neighborhood of $x$.
By lemma \ref{lem:coord} there are local coordinates $(x,t)$ near $x$ such that
$\kO=\ensemble{(x,t)}{t=0}$ and $\gth.v=\gh.v$ with
$\gh=\sum_{i=1}^s(\frac{\gl_i}2+1)t_iD_{t_i}$. By corollary \ref{cor:support}, there is
a polynomial $b_x$ with roots lower or equal to
$-Tr(\gh)=-(\gl_\kp(x)+\dim\kp)/2$ such that $b_x(\gh).v=b_x(\gth).v$ vanishes
on a neighborhood of $x$.
\end{proof}

In the diagonal case, we get:
\begin{cor}\label{cor:diag}
Let $\kg$ be a semisimple Lie algebra and $F$ be a graduate and finite
codimensional ideal of $S(\kg)^G$.

The roots of the $b$-function of $\MFf$ at $\{0\}$ are lower or equal
to $-(\rank\kg+\dim\kg)/2$ and the roots of the $b$-function of $\MF$ at
$\{0\}$ are greater or equal to $(\rank\kg-\dim\kg)/2$.
\end{cor}

\begin{proof}
In the case of a semisimple Lie group $G$ acting on its Lie
algebra $\kg$ we have $\dim\kg=\sum_{i=1}^r(\gl_i+1)$ where $r$ is
the codimension of the orbit of $x$ \cite[Ch 5.6.]{VARAD2} hence,
by definition, $\gl_\kp(x)=r$ . The minimum of $\gl_\kp(x)$ is
thus the rank of $\kg$. The result on $\MFf$ gives the
corresponding result on $\MF$ by proposition \ref{prop:fourier}.
\end{proof}

In the diagonal case, the numbers $\gl_\kp$ (defined in section
\ref{sec:Sym}) and $\tilde\gl_\kp$ are equal. In the general case,
we will use the fact that a non distinguished nilpotent point
commutes with semisimple points to get a better result than
proposition \ref{prop:grad}. Before doing this we have to study
the module $\MF$ in a neighborhood of a semisimple point.

The main property of $\MF$ is to be constant along the orbits of $K$
but also on the center of $\kg$.  Let us recall that if $X$ is a
complex manifold equal to a product $X=Y\x Z$ and $\CN$ is a coherent
$\DZ$-module, the external product of $\OY$ by $\CN$ is by definition
the coherent $\DX$-module
$$\OY\widehat\ox\CN=\DX\ox_{(q^{-1}\DY\ox_\C\ p^{-1}\DZ)}
(q^{-1}\OY\ox_\C\ p^{-1}\CN)$$
where $p:Y\x Z\to Z$ and $q:Y\x Z\to Y$ are the canonical
projections. Remark that $\OY\widehat\ox\CN$ is equal to the inverse image
$p^*\CN$. We say that a $\DX$-module is constant along $Y$ if it is
isomorphic to a module of this form.

We assume now that $\kg$ is reductive. Let $\kc$ be the center of
$\kg$ and set $\tkg=[\kg,\kg]$, $\tkp=\kp\cap\tkg$, $\tkk=\kk\cap\tkg$
and $\kc_\kp=\kp\cap\kc$. Then $(\tkg,\tkk)$ is a symmetric pair with
$\tkg$ semisimple and $\tkg=\tkk\oplus\tkp$. We have
$\kg=\kc\oplus\tkg$, $\kp=\kc_\kp\oplus\tkp$ and this decomposition is
compatible with the stratifications of $\kp$ and $\tkp$ defined in
\S\ref{sec:stratificationbis}. The action of $K$ on $\kc_\kp$ being
trivial, $S(\kp)^K=S(\kc_\kp)\ox S(\tkp)^K$. This defines a graduate
morphism $\gd_c:S(\kp)^K\to S(\tkp)^K$ by restriction and, $\gd_c$
being surjective, $F_c=\gd_c(F)$ is an ideal of finite codimension of
$S(\tkp)^K$. Let $\CM_{F_c}=\CaD_{\tkp}/\CI_c$ where $\CI_c$ is the
ideal of $\CaD_{\tkp}$ generated by $\gt_{\tkp}(\kk)$ and $F_c$. We
proved in \cite[lemma 2.2.3.]{LBY} that $\MF$ is isomorphic to
$\CO_{\kp_0}\widehat\ox(\MF)_{\tkp}$ and that $(\MF)_{\tkp}$ (the
inverse image of $\MF$ on $\tkp$) is a quotient of a power of
$\CM_{F_c}$. Concerning the $b$-functions we have:

\begin{lem}\label{lem:ss} Let $\gS$ be a submanifold of $\tkp$ and
$\gh$ be a vector field definite positive with respect to $\gS$. Let $b$
be a $b(\gh)$-function for $\CM_{F_c}$ along $\gS$.

1) $b$ is a $b(\gh)$-function for $\MF$.

2) Let $\gth_0$ be the Euler vector field of $\kc_\kp$. There exists
some $N\in\N$ such that $b(T)b(T-1)\dots b(T-N)$ is a
$b(\gh+\gth_0)$-function for $\MF$.
\end{lem}

\begin{proof}
Assume that $\kp=\kp_0\oplus\kp_1$, the action of $K$ on $\kp_0$ being
trivial and denote by  $\gd_1:S(\kp)^K\to S(\kp_1)^K$ the restriction
morphism, $F_1=\gd_1(F)$,  $\CI_1$ and $\CM_{F_1}$ the corresponding
modules. We will prove the lemma in this more general situation.

We may assume that $\kp_0=\C$ and choose linear coordinates
$(x,t)$ of $\kp$ such that $\kp_0=\ensemble{(x,t)\in\kp}{x=0}$.
Then we identify $\CaD_{\kp_1}$ to the subsheaf of $\Dp$ of
differential operators independent of $(t,D_t)$ and
$\gt_{\kp_1}(\kk)$ corresponds to $\gt_{\kp}(\kk)$. As $F_1$ is
identified to a subset of $F$, $\CI_1$ is a subsheaf of $\CI$. For
this immersion, the $V^\gh\CaD_{\kp_1}$-filtration is compatible
with both the $V^\gh\Dp$ and the $V^{\gh+\gth_0}\Dp$-filtrations.

If $b$ is a $b(\gh)$-function for $\CM_{F_1}$, this means that $\CI_1$
contains an operator $b(\gh)+Q$ with $Q\in V_{-1}^\gh\CaD_{\kp_1}$,
this gives immediately a $b(\gh)$-function for $\CM_F$.

The action of $K$ is trivial on $\kp_0$ hence $S(\kp)^K$ contains
$S(\kp_0)$. As $F$ is finite codimensional in $S(\kp)^K$ it contains a
polynomial in the dual variable of $t$ that is a polynomial in the
differential operator $D_{t}$.

Denote $\gth_0=tD_t$ the Euler vector field of $\kp_0$. Let
$A(D_t)=D_t^N+a_1D_t^{N-1}+\dots+a_N$ be the polynomial in $F$
hence $t^NA(D_t)$ is in $F$ and $t^NA(D_t)=t^ND_t^N+tR(t,tD_t)=
\gth_0(\gth_0-1)\dots(\gth_0-N+1)+tR(t,tD_t)$. We have
$b(\gh+\gth_0)=b(\gh)+g(\gh,\gth_0)\gth_0$, hence
\begin{align*}
b(\gh+\gth_0)b(\gh+\gth_0-1)\dots b(\gh+\gth_0-N+1)&=\\
c_N(\gh,\gth_0)b(\gh)+g_N(\gh,\gth_0)\gth_0&(\gth_0-1)\dots(\gth_0-N+1)
\end{align*}

This means that $b(\gh+\gth_0)b(\gh+\gth_0-1)\dots
b(\gh+\gth_0-N+1)$ is in the graduate of $\CI$ for the
$V^{\gh+\gth_0}$-filtration and shows the second part of the
lemma.
\end{proof}

In the next proposition, we assume again that $\kg$ is semisimple.

Let $s$ be a non-zero semisimple element of $\kp$. Then
$\kp=\kp^s\oplus[\kk,s]$ and $\kg^s=\kk^s\oplus\kp^s$ defines a
symmetric pair.  Set $\tkg=[\kg^s,\kg^s]$, $\tkp=\kp^s\cap\tkg$,
$\tkk=\kk^s\cap\tkg$. Then $(\tkg,\tkk)$ is a symmetric pair with
$\tkg$ semisimple and $\tkg=\tkk\oplus\tkp$. If $\kc$ is the center of
$\kg^s$, $\kp^s=\tkp\oplus(\kp^s\cap\kc)$.

Let $\ka$ be a Cartan subspace of $\kp$ containing $s$,
$\gF=\gF(\kg,\ka)$ the root space,
$P=\ensemble{\ga\in\gF}{\ga(s)=0}$. Then $\kp^s=\ka\oplus\kp_P$,
$\ka=\ka_P\oplus\kapo$ and the stratum of $s$ is
$$S_{(P,\{0\})}=\bigcup_{x\in\kapp}K.x=K.\kapp$$

\begin{prop}\label{prop:semi}
There are local coordinates
$(x_1,\dots,x_l,y_1,\dots,y_p,t_1,\dots,t_q)$ of $\kp$ such that:

(i)\ $\kp^s=\ensemble{(x,y,t)}{x=0}$,
$\tkp_s=\ensemble{(x,y,t)}{x=0,y=0}$, $s=(0,y_0,0)$ with
$y_0\ne0$.

(ii)If $z$ is a point of $\kp$ close to $s$ whose semisimple part
is $s$, the stratum of $z$ is equal in a neighborhood of $s$ to
the set $S_{(P,\kO)}=\ensemble{(x,y,t)}{t\in\kO}$ where $\kO$ is
the orbit of the nilpotent part of $z$ in $\tkp_s$. The stratum of
$s$ is $S_{(P,\{0\})}=\ensemble{(x,y,t)}{t=0}$.

(iii) The vector fields $D_{x_1},\dots,D_{x_l}$ are in the ideal of $\Dp$
generated by $\gt(\kk)$.

(iv) The Euler vector field $\gth$ of $\kp$ is equal to
$$\gth=\sum_{i=1}^py_iD_{y_i}+\sum_{j=1}^qt_jD_{t_j}$$
\end{prop}

\begin{proof}
Let $V$ be a subspace of $\kk$ such that $\kk=V\oplus\kk^s$. Let
$(u_1,\dots,u_l)$ be a basis of $V$, $(v_1,\dots,v_p)$ be a basis
of $\kapo$, $(w_1,\dots,w_q)$ be a basis of $\tkp_s$. Let $y_0\ne
0$ be the coordinate of $s$ in the basis $v$. We have
$\kp=[\kk,s]\oplus\kapo\oplus\tkp_s$ hence the map:
$$F(x,y,t)= \exp(x_1u_1)\dots\exp(x_lu_l).(\sum y_iv_i + \sum t_jw_j)$$
defines an isomorphism from a neighborhood of $(0,y_0,0)$ to a
neighborhood of $s$ in $\kp$ hence defines local coordinates of
$\kp$. These coordinates satisfy the condition (i) by definition.
As $\kk=V\oplus\kk^s$, in a neighborhood of $s$ the orbits of $K$
are of the form $\ensemble{(x,y,t)}{y=c, t\in K^sd}$ for some
$c\in\kapo$ and $d\in\tkp_s$ which shows (ii).Next, (iii) is
satisfied from \cite[lemma 3.7.]{SEKI}.

Let us calculate the Euler vector field $\gth$ of $\kp$ in these
coordinates. By definition, $\gth$ acts as $\gth(f)(z)=\frac
d{d\ga}f(e^{\ga}z)|_{\ga=0}$ in linear coordinates $z$ of $\kp$, hence in
coordinates $(x,y,t)$:
$$\gth(f)(x,y,t)=\frac d{d\ga}f(F^{-1}(e^{\ga}F(x,y,t)))|_{\ga=0}$$
The $K$-action commutes with scalar multiplication hence
$e^{\ga}F(x,y,t)=F(x,e^{\ga}y,e^{\ga}t)$ and thus
$\gth=\sum_{i=1}^py_iD_{y_i}+\sum_{j=1}^qt_jD_{t_j}$.
\end{proof}

Let $\gd_s$ be the restriction map $S(\kp)^K\to S(\kp^s)^{K^s}$ which
is graduate. If $F$ is an ideal of finite
codimension of $S(\kp)^K$, the set of points of $\kp^*$ defined by $F$
is a finite union of orbits of $\kp^*$ hence its intersection with
$(\kp^s)^*$ is also a finite union of orbits hence $\gd_s(F)$ is an ideal of finite
codimension of $S(\kp^s)^{K^s}$. Let $\CI_s$ be the left ideal of
$\Dps$ generated by $\gd_s(F)$ and $\gt(\kk^s)$ and
$\CM_s=\Dps/\CI_s$. In a neighborhood of $s$ we identify $\kp$ to the
product of the orbit $K.s$ by $\kp^s$, then we have:

\begin{prop}\label{prop:semib}
In a neighborhood of $s$, the module $\MF$ is isomorphic to
$\CO_{Ks}\widehat\ox\CN$ where $\CN$ is a quotient of $\CM_s$.

Let $\gS_0$ be a stratum of $\kp^s$ and $\gh$ be a vector field
on $\kp^s$ which is definite positive with respect to $\gS_0$. Let $b$ be
a polynomial which is a $b(\gh)$-function for $\CM_s$. Then
$\gS=K.s\x\gS_0$ is a stratum of $\kp$ in a neighborhood of $s$,
$\gh$ is definite positive with respect to $\gS$ and $b$ is a
$b(\gh)$-function for $\MF$.
\end{prop}

\begin{proof}
We use the local coordinates $(x,y,t)$ of lemma \ref{prop:semi}. By (iii) of
this lemma, the vector fields $D_{x_1},\dots,D_{x_l}$ are in the ideal $\CI_F$
hence $\MF$ is isomorphic to $\CO_{Ks}\widehat\ox\CN$ for some coherent
$\Dps$-module $\CN=\Dps/\CJ$.

Let $a\in\kk^s$, $f\in S(\kp^*)$, $f_s$ the restriction of $f$ to
$\kp^s$ and $\gt_{\kp^s}(a)$ the vector field associated to $a$ by the
action of $K^s$ on $\kp^s$. By definition, if $u\in\kp^s$:
$$(\gt_{\kp^s}(a)f_s)(u)=\frac d{ds}f\left(\exp(-ta).u\right)|_{s=0}$$
is equal to the restriction of $\gt_{\kp}(a)f$ to $\kp^s$ hence
$\gt_{\kp}(a)=\gt_{\kp^s}(a)+w$ where $w$ is a vector field on $\kp$ vanishing
on $\kp^s$. As the ideal $\CI_F$ contains $D_{x_1},\dots,D_{x_l}$ and
$\gt_{\kp}(a)$ it contains $\gt_{\kp^s}(a)$. This means that $\CJ$ contains
$\gt_{\kp^s}(\kk^s)$.

On the other hand, let $P\in F$, as the coordinates $(y,t)$ are
linear coordinates of $\kp^s$, the value of $P$ on a function of
$t$ is the restriction of $P$ to $S(\kp^s)^{K^s}$. Hence $\CJ$
contains $\gd_s(F)$ and $\CN$ is a quotient of $\CM_s$.

The second part of the proposition is clear for $\CO_{Ks}\widehat\ox\CM_s$
hence for $\CM$.
\end{proof}

Recall that $\gl_\kp$ is the minimum of $\gl_\kp(x)$ for all
distinguished nilpotents $x$ and if $s$ is a semisimple element of
$\kp$,  $\gl_{\kp^s}$ is defined in the same way with
$\kp^s=\ensemble{x\in\kp}{[x,s]=0}$. We defined also $\gm_\kp$ as the
minimum over all semisimple elements $s\in\kp$
of $(\gl_{\kp^s}-\redim\kp^s)/2$.

\begin{prop}\label{prop:precis}
The roots of the $b$-function of $\MFf$ at $\{0\}$ are lower or equal
to $-\gm_\kp-\dim\kp$ and the roots of the $b$-function of $\MF$ at
$\{0\}$ are greater or equal to $\gm_\kp$.
\end{prop}

\begin{proof}
We will prove the proposition by induction, assuming that the result
has been proved for all the symmetric sub-pairs of $(\kg,\kk)$. We keep the
notations of the proof of proposition \ref{prop:grad} and make the
same proof except for non distinguished nilpotent points.

Let $u$ be the canonical generator of $\MF$ and $v$ be the canonical generator
of $\MFf$. Let $b_0$ be a polynomial such that $b_0(\gth^*)v$ is supported by
a nilpotent orbit $\kO$ in a neighborhood of one of its points $x$. We assume
that the roots of $b_0$ are lower or equal to $-\gm_\kp-\dim\kp$ and that $x$
is not distinguished (otherwise we use the proof of \ref{prop:grad}).

By definition of non distinguished points, there exists some
$s\in\kp$ which is semisimple and such that $[x,s]=0$. Then by
proposition \ref{prop:semib} and the induction hypothesis, there
is a polynomial $b_1$ with roots greater or equal to $\gm_\kp$
such that $b_1(\gth)u$ vanishes at $s$. Remark that all semisimple
points of $\tkp$ are semisimple in $\kp$ hence
$\gm_{\tkp}\ge\gm_\kp$.

Proposition \ref{prop:fourier} shows that the characteristic variety of the
module generated by $b_1(-\gth^*-\dim\kp)v$ does not contain the point $(x,s)$.
The module generated by $b_0(-\gth^*-\dim\kp)b_1(-\gth^*-\dim\kp)v$ is
supported by $\kO$ and $x$ is a smooth point of $\kO$, hence if $x$ is in the
support of the module, the conormal bundle to $\kO$ at $x$ is contained in the
characteristic variety. But $(x,s)$ is a point of this conormal bundle which
does not belong to the characteristic variety, hence
$b_0(-\gth^*-\dim\kp)b_1(-\gth^*-\dim\kp)v$ vanishes at $x$.
\end{proof}

Now, we do not assume any more that $F$ is graduate, then

\begin{cor}\label{cor:nongrad}
The roots of the $b$-function of  $\MF$ at $\{0\}$ are greater or
equal to  $\gm_\kp$.
\end{cor}

\begin{proof}Let $F'$ be the graduate of $F$, then by proposition
\ref{prop:precis}, we have an equality $b(\gth)=\sum A_i(x,D_x)
u_i(x,D_x) + B_j(x,D_x) Q_j(D_x)$ where $A_i$ and $B_j$ are
differential operators of $\Dp$, $u_i$ are vector fields of $\gt(\kk)$
and $Q_j\in F'$.

We have $[\gth,u]=0$ for any $u$ in $\gt(\kk)$ hence $u$ is of degree
$0$ for the graduation associated to the $V$-filtration along
$\{0\}$. On the other hand, if $Q\in F'$ is homogeneous of
degree $k$ as a polynomial, we have $[Q(D_x),\gth]=kQ(D_x)$ that is
$Q$ is of degree $k$ for the $V$-filtration. Decomposing $A_i$ and
$B_j$ in homogeneous parts, we may rewrite $b(\gth)=\sum \tilde
A_i(x,D_x) u_i(x,D_x) + \tilde B_j(x,D_x) Q_j(D_x)$ with $\tilde
A_i(x,D_x) u_i(x,D_x)$ and $\tilde B_j(x,D_x) Q_j(D_x)$ homogeneous of
degree $0$ for the $V$-filtration.

Now if $Q_j\in F'$, there exists $P_j\in F$ such that $P_j=Q_j+R_j$
with $R_j$ of degree lower than the degree of $Q_j$ hence
$$b(\gth)+ \sum \tilde B_j(x,D_x) R_j(D_x)=\sum \tilde
A_i(x,D_x) u_i(x,D_x) + \tilde B_j(x,D_x) P_j(D_x)$$ which means that
$b(\gth)+ \sum \tilde B_j(x,D_x) R_j(D_x)$ is a $b$-function.
\end{proof}

Let $b$ be the $b$-function of $\MF$ at $\{0\}$ and for each nilpotent
orbit $\kO$ of $\kp$ let $\gh_\kO$ defined by lemma \ref{lem:coord}.

\begin{prop}\label{prop:nilp}
For each nilpotent orbit $\kO$ of $\kp$, $b(\gh_\kO)$ is a
quasi-$b$-function for $\MF$. If $F$ is a graded ideal, this
quasi-$b$-function is monodromic.
\end{prop}

\begin{proof}
By the hypothesis, $b$ is the $b$-function of $\MF$ at $\{0\}$ hence
the ideal $\CI_F$ contains an equation $b(\gth)+R$ where $R$ is a
differential operator of order $-1$ for the $V$-filtration at
$\{0\}$. By lemma \ref{lem:coord} $\gh_\kO=\gth$ in $\CI_F$, hence
$b(\gh_\kO)+R$ is also in $\CI_F$. For simplicity, we will write $\gh$
for $\gh_\kO$ in this proof.

As $R$ is of order $-1$ for the $V$-filtration, we can write it as a
series $R=\sum_{k\le-1}R_k$ with $[R_k,\gth]=kR_k$. Let
$R_k(x,t,D_x,D_t)=R_k^0(x,t,D_t)+\sum R_k^i(x,t,D_x,D_t)D_{x_i}$, we
have  $[R_k,\gth]=[R_k^0,\gth]$ modulo $D_x$ hence
$[R_k^0,\gth]=kR_k^0$. As $D_{x_1},\dots,D_{x_{n-r}}$ are in the ideal
$\CI_F$ by lemma \ref{lem:coord}, we may replace $R$ by $\sum_{k\le-1}R_k^0$ and assume from now
on that $R$ is independent of $D_x$.

We decompose now each $R_k$ as a series
$R_k(x,t,D_t)=\sum_jR_{kj}(x,t,D_t)$ where each $R_{kj}$ is
homogeneous of degree $j$ for $\gh$, that is
$[R_{kj},\gh]=jR_{kj}$. By uniqueness of the decomposition, we
have $[R_{kj},\gth]=kR_{kj}$ hence
$[R_{kj},D_{x_{n-r}}]=2(k-j)R_{kj}$, that is
$R_{kj}=R_{kj}^1(x',t,D_t)e^{2(j-k)x_{n-r}}$ with
$x'=(x_1,\dots,x_{n-r-1})$. Finally $R$ is equal to a convergent
series:
$$R(x,t,D_t)=\sum_{k\le-1,j\le j_0}R_{kj}(x',t,D_t)e^{2(j-k)x_{n-r}}$$
where $R_{kj}$ is homogeneous of degree $j$ for $\gh$. The ideal
$\CI_F$ contains the operators $D_{x_i}$ for $i=1,\dots,n-r$ hence is
generated by these $D_{x_i}$ and by a finite number
$Q_1(t,D_t),\dots,Q_N(t,D_t)$ of differential operators independent of
$(x,D_x)$ and thus we have:
$$b(\gh)+R(x,t,D_t)=\sum_{i=1}^NA_i(x,t,D_t)Q_i(t,D_t)$$
Therefore, the operator $b(\gh)+R(0,x_{n-r},t,D_t)$ obtained by making
$x'=0$ is still in $\CI_F$. In the same way, the operator obtained by
integration on the path $x_{n-r}\in[0,2i\pi]$ is still in $\CI_F$. But
$\int_{[0,2i\pi]}e^{2(j-k)u}du = 2i\pi$ if $j=k$ and $0$ otherwise,
hence the operator
$$b(\gh)+\sum_{k\le-1}R_{kk}(0,t,D_t)$$
is an operator of $\CI_F$. By construction,
$\sum_{k\le-1}R_{kk}(0,t,D_t)$ is a differential operator of order
$-1$ for the $V^\gh$-filtration, hence $b(\gh)+\sum_{k\le-1}R_{kk}(0,t,D_t)$
is a quasi-$b$-function.

If $F$ is a graded ideal, we have $R=0$ from corollary \ref{prop:grad}
hence $b(\gh)$ is a monodromic quasi-$b$-function.
\end{proof}

\begin{proof}[Proof of theorem \ref{thm:mainsym}]
We defined in \S\ref{sec:stratificationbis} a finite stratification of
$\kp$. We will define a vector field $\gh_\gS$ definite positive with
respect to $\gS$ with trace equal to $t_\gS$ and show that $\MF$
admits a $b(\gh_\gS)$-function whose roots are greater or equal to
$\gm_\gS$.

Assume first that the theorem has been proved when $\kg$ is semi-simple.
If $\kg$ is not semi-simple, we set as in remark \ref{rem:prod}
$\tkg=[\kg,\kg]$, $\tkp=\kp\cap\tkg$, $\tkk=\kk\cap\tkg$ and
$\kc_\kp=\kp\cap\kc$ where $\kc$ is the center of $\kg$. The strata of
$\kp$ are equal to the direct sum of $\kc_\kp$ and the strata
of $\tkp$. Let $\gS$ be a stratum of $\tkp$, we associate to
$\kc_\kp\oplus \gS$ the same vector field $\gh_\gS$ and lemma
\ref{lem:ss} gives the result.

So, we may assume now that $\kg$ is semisimple, take $x\in\kp$, and
prove the result for the stratum $\gS$ of $x$. Take first $x=0$. Then
$\gS=\{0\}$, $\gh_\gS=\gth$ the Euler vector field of $\kp$ with trace
$\dim\kp$ and corollary \ref{cor:nongrad} shows that $\MF$ admits a
$b$-function $b_0$ whose roots are greater or equal to $\gm_\kp$.

Assume now that $x$ is a nilpotent point of $\kp$. Then $\gS$ is the
orbit of $x$, $\gh_\gS$ is the vector field defined by lemma
\ref{lem:coord} whose trace is $t_\gS$ and proposition
\ref{prop:nilp} shows that $b_0$ is a  $b(\gh_\gS)$-function for
$\MF$.

Consider now a non-nilpotent point $x$ with Jordan decomposition
$x=s+n$. We may assume by induction on the dimension of $\kp$ that
the theorem has been proved for the pair $\kg^s=\kk^s\oplus\kp^s$.
As in the proof of  lemma \ref{prop:stratp}, we may assume that
$x$ is arbitrarily close to $s$. Then the result follow from
proposition \ref{prop:semib}.

To end the proof of the theorem, we remark that if $F$ is graduate,
all $b$-functions are monodromic, this shows that the singular support
$\kp\setminus\kprs$ of $\MF$ is conic relatively to all vector fields
$\gh_\gS$ by remark \ref{geom-monodrom}. As they do not depend on $F$
the result is still true if $F$ is not graduate.
\end{proof}

The other results of sections \ref{sec:HC} and \ref{sec:Sym} are
direct consequences of theorem \ref{thm:mainsym}:

\begin{proof}[Proof of corollary \ref{cor:premain}]
Theorem \ref{thm:mainsym} shows that $\MF$ is conic-tame if for
any stratum $\gS$ we have $\gm_\gS+t_\gS>0$, that is for
$x=s+n\in\gS$ if $\gm_{\kp^s}+(\gl_{\kp^s}(n)+\redim\kp^s)/2>0$.
If $\lpx>0$ for any  $x$, this is true by definition of
$\gm_{\kp^s}$.
\end{proof}

\begin{proof}[Proof of corollary \ref{cor:wpremain}]
Consider now a nilpotent point $x$ of $\kp$, $\kO$ its orbit. If
$T^*\kp$ is identified to $\kp\x\kp$, the conormal bundle to $\kO$
is $\ensemble{(x,y)\in \kO\x\kp}{[x,y]=0}$ and the characteristic
variety of $\CM$ is contained in
$\ensemble{(x,y)\in\kp\x\kp}{[x,y]=0, y\in \kN(\kp)}$. If $x$ is
not distinguished, there exists some semi-simple $y$ such that
$[x,y]=0$ and if $\gp$ is the projection
$T^*\kp\simeq\kp\x\kp\to\kp$, $\gp^{-1}(x)\cap\CHM$ is strictly
contained in the conormal bundle to $\kO$. This is true for all
points $x'$ of $\kO$ and thus $\kO$ satisfies the condition of
definition \ref{def:wtame}(ii). If $x=s+n$ is the Jordan
decomposition of $x\in\kp$ and if $n$ is not distinguished in
$\kp^s$, the same condition is still true for the stratum of $x$.
If we assume only that $\gl_{\kp^s}>0$ for any $s\in\kp$
semisimple, we get the inequality $\gm_\gS+t_\gS>0$ for all
$x=s+n$ such that $n$ is distinguished hence $\CM$ is weakly tame.
\end{proof}

\begin{proof}[Proof of theorem \ref{thm:main}]
In the diagonal case, $\gl_{\kp^s}$ is always strictly positive
which shows the theorem.
\end{proof}

 Then corollaries \ref{cor:main},
\ref{cor:mainsym} and \ref{cor:secsym} are deduced from the
results of section \ref{sec:tameres}.

\begin{proof}[Proof of proposition \ref{prop:prec}]
In the diagonal case, if $\gS$ is a nilpotent orbit and $x\in\gS$,
the trace of $\gh_\gS$ is $t_\gS=(\gl_\kg(x)+\dim\kg)/2\ge
(\rank\kg+\dim\kg)/2$ while corollary \ref{cor:diag} shows that
the roots of the $b$-function of $\gS$ are greater or equal to
$(\rank\kg-\dim\kg)/2$. So, the roots of the $b$-function are
greater or equal to $-t_\gS/\gd$ if
$\gd=(\dim\kg+\rank\kg)/(\dim\kg-\rank\kg)$. If $\gd(\kg)$ is the
minimum of this value over all semi-simple subalgebras
$[\kg^s,\kg^s]$ for $s$ semi-simple, the roots of the $b$-function
of $\gS$ will be greater or equal to $-t_\gS/\gd(\kg)$ for all
strata $\gS$ and definition \ref{def:tame} will be satisfied.
\end{proof}

\providecommand{\MR}{\relax\ifhmode\unskip\space\fi MR }
\providecommand{\MRhref}[2]{%
  \href{http://www.ams.org/mathscinet-getitem?mr=#1}{#2}
}
\providecommand{\href}[2]{#2}

\enddocument

\end